\documentclass{article}
\usepackage{graphicx}
\usepackage{arxiv}

\usepackage[utf8]{inputenc} % allow utf-8 input
\usepackage[T1]{fontenc}    % use 8-bit T1 fonts
\usepackage{url}            % simple URL typesetting
\usepackage{booktabs}       % professional-quality tables
\usepackage{amsfonts}       % blackboard math symbols
\usepackage{nicefrac}       % compact symbols for 1/2, etc.
\usepackage[dvipsnames]{xcolor}
\usepackage{physics}
\usepackage{graphicx} % Required for inserting images
\graphicspath{ {./images} }
\usepackage{amsmath}
\usepackage{amssymb}
\usepackage{amsthm}
\usepackage{mathtools}
\usepackage{float}
\usepackage{todonotes}
\usepackage{array}
\usepackage{enumerate}
\usepackage{multicol}
\usepackage{multirow}
\usepackage{mathrsfs}
\usepackage{siunitx}
\usepackage{tabularx}
\usepackage{subcaption}
\usepackage[bookmarks, hidelinks]{hyperref}

\usepackage{bbm}
\usepackage{authblk}

\usepackage{subcaption}
\usepackage{hyperref}
\usepackage{booktabs}
\usepackage{cleveref}

\usepackage{adjustbox}
\usepackage{amsopn}

\def\fp{\delta f}

\def\up{\delta u}

\def\Diff{D}
\def\Di{\mathcal{D}^s}

\AtBeginDocument{\RenewCommandCopy\qty\SI}

\makeindex  

\theoremstyle{definition}
\newtheorem{definition}{Definition}[section]

\newtheorem{remark}{Remark}

\newcommand{\inp}[2]{\left\langle #1,\, #2 \right\rangle}

\renewcommand{\vec}[1]{\mathbf{#1}}

\DeclareMathOperator{\diag}{diag}

\newcommand{\lJump}{[\![}
\newcommand{\rJump}{]\!]}

\renewcommand{\vec}[1]{\mathbf{#1}}

\newcommand{\fn}[2]{\mathinner{#1\mathopen{\left(#2\right)}}}

\title{Local linear stability of dual-pairing summation-by-parts methods for nonlinear conservation laws}

\author[1]{Dougal Stewart}
\author[2]{Kenneth Duru}

\affil[1]{School of Mathematics, Monash University, Australia}
\affil[2]{Department of Mathematical Sciences, University of Texas at El Paso, USA}

\date{}

\begin{document}

\maketitle

 \begin{abstract}
A recent study by Gassner et al. [J. Sci. Comput. 90:79 (2022)] demonstrates that local energy stability---that is, ensuring the asymptotic numerical growth rate does not exceed the continuous growth rate---is crucial for achieving accurate numerical simulations of nonlinear conservation laws. While nonlinear entropy stability is necessary for numerical stability (i.e., ensuring the boundedness of nonlinear numerical  solutions), local energy stability is essential to prevent unresolved high-frequency wave modes from dominating the simulation.
Currently, it remains an open question whether high-order numerical methods for nonlinear conservation laws can be simultaneously entropy-stable and locally energy-stable. In this work, we examine the local energy-stability properties of recently developed entropy-stable, high-order accurate dual-pairing (DP) SBP methods, as introduced by Duru et al. [arXiv: 2411.06629], for nonlinear conservation laws. Our analysis indicates that the entropy-stable volume upwind filter inherent in these methods can ensure local energy stability.
This approach offers a novel numerical strategy for designing reliable high-order methods for nonlinear conservation laws that are provably entropy-stable and locally energy-stable. The theoretical findings are supported by numerical experiments involving the inviscid Burgers’ equation and nonlinear shallow water equations, in one and two space (2D) dimensions. Furthermore, we present accurate numerical simulations of 2D barotropic shear instability, with fully developed turbulence, demonstrating the efficiency of the DP SBP method in resolving turbulent scales.
 \end{abstract}

\section{Introduction}
Enhancing the robustness of high-order accurate numerical methods for nonlinear conservation laws remains a fundamental challenge. Historically, many studies have promoted schemes with provable linear stability, leveraging assumptions of solution smoothness to approximate nonlinear problems \cite{H-O-Kreiss1970,kreiss2004initial,Gustafsson-Sundstrom1978,Oliger-Sundstrom1978}. However, linear stability alone often proves insufficient for capturing complex nonlinear phenomena such as turbulence and shock formation, where nonlinear effects dominate. Empirical evidence from numerous numerical simulations \cite{Tadmor2003entropystable,gassner2022stability,FISHER2013353,ChanRanochaHendrik_et_al2022} demonstrates that ensuring nonlinear numerical stability is crucial for robustness of high-order numerical methods for nonlinear conservation laws.

For entropy-stable partial differential equations (PDEs), constructing provably nonlinearly stable numerical methods can be achieved by replicating continuous entropy stability properties at the discrete level. A key tool in this endeavor is the summation-by-parts (SBP) principle \cite{kreiss1974finite,BStrand1994,fernandez2014review,svard2014review,gassner2013skew,taylor2010compatible}, which possesses a mimetic structure facilitating the design of stable semi-discrete schemes, provided boundary conditions are carefully managed. At the continuous level, entropy stability analysis uses integration by parts, along with the chain rule and product rule. Although SBP operators discretely mimic integration by parts, replicating the chain and product rules remains challenging. To mirror continuous analysis, nonlinear conservation laws are often reformulated into skew-symmetric or entropy-conserving split-forms, circumventing direct application of these rules at the discrete level and enabling stability proofs, see \cite{Sjögreen2019,GERRITSEN1996245,FISHER2013353,STIERNSTROM2021110100,SHU1988439,ricardo2024entropy,ricardo2023entropy,ChanRanochaHendrik_et_al2022,gassner2016skewsym_swe,Nordstrom_2022_skewsym_comp}.

Once an appropriate entropy functional and an entropy-conserving reformulation of the PDE system are established, the system can be discretized using SBP operators to produce entropy-conserving schemes. Although this approach guarantees nonlinear stability in principle, traditional SBP finite difference (FD) methods \cite{kreiss1974finite,BStrand1994,fernandez2014review,svard2014review} and collocated discontinuous Galerkin spectral element methods (DGSEM) \cite{gassner2013skew,taylor2010compatible,fluxo_SplitForm} are often prone to generating spurious unresolved wave modes. These artifacts can degrade solution accuracy or even cause numerical instabilities. Recent research \cite{gassner2022stability} emphasizes that local linear stability of high-order methods is vital for reliable numerical simulations of conservation law. Unfortunately, many existing high-order entropy-stable schemes—whether SBP FD, DGSEM, or spectral collocation methods—implemented in skew-symmetric split form, typically lack this local linear stability.

The dual-pairing (DP) SBP framework \cite{DovgilovichSofronov2015,MATTSSON2017upwindsbp,williams2024drp} has recently been developed to improve the accuracy of high-order methods for wave propagation and conservation laws. This approach has been extended to discontinuous Galerkin methods \cite{GLAUBITZ2025113841}, enabling more flexible and robust treatment of complex geometries. DP-SBP operators consist of a dual pair of high-order discrete operators that satisfy the SBP property, ensuring stability and consistency. These operators include additional degrees of freedom, which can be optimized to reduce dispersion errors \cite{williams2024drp}.

The DP-SBP method for nonlinear conservation laws was recently introduced in \cite{duru2024dualpairingsummationbypartsfinitedifference}. It integrates three key components: 1) the DP SBP DG/FD operators, 2) skew-symmetric reformulation, and 3) upwind flux splitting. This method is conservative, entropy consistent (correctly capturing entropy loss at shocks),  provably entropy-stable and naturally incorporates intrinsic high-order accurate filters  at shocks and discontinuities. In this work, we analyze the local linear stability properties of the DP-SBP method \cite{duru2024dualpairingsummationbypartsfinitedifference} for nonlinear conservation laws. Our analysis indicates that the entropy-stable volume upwind filter embedded within this approach can ensure local linear stability.

This framework offers a novel numerical strategy for designing reliable high-order methods for nonlinear conservation laws that are both provably entropy-stable and locally linear stable. The theoretical insights are supported by numerical experiments involving the inviscid Burgers’ equation and nonlinear shallow water equations (SWEs) in one and two spatial dimensions. Additionally, we present high-fidelity simulations of 2D barotropic shear instability with fully developed turbulence, demonstrating the DP SBP method’s effectiveness in resolving turbulent scales.

The paper is organized as follows. Section\ref{sec:continuous-and-discrete-analysis} introduces a scalar model problem and discusses entropy stability and local linear stability concepts. Section\ref{sec:sbp-and-scalar} details the DP-SBP framework, presents the local linear stability analysis, and derives optimal upwind parameters for Burgers’ equation. Numerical experiments validating these findings are provided in Section\ref{sec:instructive-numerical-experiments}. Section\ref{sec:swes-stability} extends the analysis to systems of nonlinear SWEs in 1D and 2D, culminating in simulations of 2D shear instability and turbulence. Finally, Section~\ref{sec:conclusion} summarizes the main conclusions of our study.

\section{Entropy stability and local linear stability}\label{sec:continuous-and-discrete-analysis}
As a model problem we consider the scalar hyperbolic nonlinear conservation law in one space dimension given by
\begin{align}\label{eq:scalar-nonlinear-conservation-law}
    \partial_t u + \partial_x f(u) = 0, \quad x \in [0, L], \quad t \in [0, T],
\end{align}
where $L>0$ is the width of the spatial domain and $T>0$ is the final time. Here $u:  [0, L]\times [0, T] \to \mathbb{R}$ is the unknown solution and the continuous function $f: \mathbb{R}\to \mathbb{R}$ is the PDE flux.
We augment the conservation law \eqref{eq:scalar-nonlinear-conservation-law} with the smooth initial condition $u(x, 0) = u_0(x)$ for $x \in [0, L]$ and periodic boundary conditions $u(0, t) = u(L, t)$ for all $t \in [0, T]$.  

\subsection{Continuous analysis}
Without loss of generality we will assume that the nonlinear conservation law \eqref{eq:scalar-nonlinear-conservation-law} admits square entropy, that is \(e = u^2\) (see for e.g. the nonlinear Burgers equation discussed in the next section), such that for sufficiently smooth solutions we have
%%%%
\begin{align}\label{eq:total-entropy-evolution}
   \frac 1 2 \dv{}{t} \int_{0}^{L} e dx= \frac 1 2 \dv{}{t} \norm{u}^2
    &=  \text{boundary terms}, \quad \forall \, t \in [0, T],
\end{align}
where the weighted $L^2([0,L])$ scalar product and norm are defined by
\begin{align}\label{eq:L2-inner-product}
(u, v)_{\rho} = \int_{0}^{L} u v \rho dx, \quad \rho>0, \quad \norm{u}^2_{\rho} = (u, u)_{\rho}.
\end{align}
If $\rho \equiv 1$ we get the standard $L^2([0,L])$ scalar product and norm, we omit the subscript ${\rho}$, and we have \(\norm{u}^2 = (u, u)\).
The boundary terms in \eqref{eq:total-entropy-evolution} vanish because of the periodic boundary conditions, thus after time integration we have the conservation of total entropy
\begin{align}\label{eq:total-entropy-conservation}
\norm{u(\cdot, T)}^2 = \norm{u_0}^2.
\end{align}
The conservation of  of total entropy implies $L_2$-stability for the solutions of the nonlinear conservation law, and can be targeted  by high-order numerical methods yielding entropy conserving/stable schemes \cite{ChanRanochaHendrik_et_al2022,duru2024dualpairingsummationbypartsfinitedifference}. However, an entropy/nonlinearly stable numerical can support  high-frequency unresolved wave modes with growing amplitude which can destroy the accuracy of numerical solutions \cite{gassner2022stability,duru2024dualpairingsummationbypartsfinitedifference}.
%%%
%\subsection{Local linear stability}

In turbulence research (see, for example, \cite{Mankbadi1994}), linear stability theory is a key mechanism that explains how laminar flow can transition to nonlinear turbulent flow. This approach involves linearizing the nonlinear system of conservation laws around a chosen base-flow state. The resulting linearized equations indicate whether the base-flow is locally stable in both space and time--meaning the flow remains laminar--or if it permits growth of disturbances, implying local instability and the potential transition to turbulence through nonlinear interactions.

This concept is closely related to linear stability in the sense defined by Strang \cite{Strang1964}. It is reasonable to expect that a numerical scheme should exhibit the same local behavior as the underlying equations. For linear PDEs, it is well-known that energy-stable schemes possess this property, as they prevent non-physical growth—the growth of the numerical solution does not exceed that of the continuous problem. In contrast, Lax-stable \cite{LaxRichtmyer1956} schemes may potentially violate this principle.

%%%
We assume the existence of a smooth, periodic base-flow \(U(x,t)\) defined on \([0, L]\) that solves the nonlinear conservation law \eqref{eq:scalar-nonlinear-conservation-law}. Our primary interest lies in examining the dynamics of a sufficiently small perturbation \(\up(x,t)\) around this base-flow, where $\epsilon\in\mathbb{R}$ and \(|\epsilon\up(x,t)| \ll |U(x,t)|\) for \((x,t) \in [0, L] \times [0, T]\).
To analyze these dynamics, we linearize the nonlinear conservation law \eqref{eq:scalar-nonlinear-conservation-law} around the smooth solution \(U(x,t)\). Specifically, we consider a perturbed solution \(u = U + \epsilon \up\), with \(0 < |\epsilon| \ll 1\), and neglect terms of order \(O(\epsilon^2)\). Expanding \(f(u)\) around \(U\), we obtain:
\(f(u) \approx f(U) + \epsilon \partial_u f|_{u=U} \up,\)
which leads to the linearized conservation law governing the evolution of the perturbation \(\up(x,t)\).
\begin{equation}\label{eq:scalar-linear-conservation-law}
    \partial_t \up + \partial_x (a \up) = 0,  \quad x \in [0, L], \quad t \in [0, T],
\end{equation}
where \( a = \partial_uf|_{u =U}\) is the Jacobian of the flux, with the initial condition \(\up(x,0) = \up_0(x)\) and periodic boundary condition \(\up(0,t) = \up(L, t)\). 
%%%
%%%

%%%
%%%
To estimate the growth/decay of the $L_2$-norm of the perturbation, we multiply the linear conservation law \eqref{eq:scalar-linear-conservation-law} by \(\up\) and integrate the product over the spatial domain \([0,L]\). By applying the product rule \(\partial_x (a \up) = a \partial_x \up + \up \partial_x a\) and integration-by-parts, we have 
%%%
\begin{equation}\label{eq:total-energy-evolution}
\begin{aligned}
    \frac 1 2 \dv{}{t} \norm{\up(\cdot, t)}^2
    &= -\int \frac{1}{2} \left(\up \partial_x (a \up) + \up a \partial_x \up + \up^2 \partial_x a \right) \dd{x} \\
    &= -\int \frac{1}{2} \up^2 \partial_x a \dd{x} \underbrace{+ a\up^2|_{x = 0} - a\up^2|_{x = L}}_{ \text{boundary terms}}.
\end{aligned}
\end{equation}
%%%
Again the boundary terms will vanish because of the periodic boundary conditions, and after time integration yields
\begin{align}\label{eq:estimate-general-base-flow}
   \norm{\up(\cdot, T)}^2 \le e^{\eta_c T} \norm{\up_0}^2, \quad \eta_c = -\min_{x,t} \partial_x a(x,t). 
\end{align}
%%%
The above estimate \eqref{eq:estimate-general-base-flow} is nearly optimal for a general base-flow $U(x,t)$. However, if we have further information about the base-flow $U(x,t)$, a sharper and more optimal estimate can be derived. For example, if the base-flow is defined such that \(a(x,t) \ne 0\) does not change sign within the domain, that is $a = a_{-}< 0$ for all \((x,t) \in [0, L]\times [0, T]\) or $a = a_{+}>0$ for all \((x,t) \in [0, L]\times [0, T]\), we have
\begin{equation}\label{eq:scalar-linear-conservation-law-strict}
    \partial_t \up \pm \partial_x (|a_{\pm}| \up) = 0,  \quad x \in [0, L], \quad t \in [0, T],
\end{equation}
where  \(|a|=|a_{+}| = a_{+}>0\) and \(|a|=|a_{-}| = -a_{-}>0\). Multiplying the linear conservation law \eqref{eq:scalar-linear-conservation-law-strict} by \(\up |a|\) and integrating the product over the spatial domain \([0,L]\), we have
\begin{equation}\label{eq:total-weighted-energy-evolution}
\begin{aligned}
    \frac 1 2 \dv{}{t} \norm{\up(\cdot, t)}^2_{|a|}
    &= \int \frac{1}{2} \up^2 \partial_t |a| \dd{x} \underbrace{\pm |a_{\pm}|^2\up^2|_{x = 0} \mp |a_{\pm}|^2\up^2|_{x = L}}_{ \text{boundary terms}}.
\end{aligned}
\end{equation}
As before, the boundary terms in \eqref{eq:total-weighted-energy-evolution} will vanish because of the periodic boundary conditions, giving
\begin{align}\label{eq:estimate-general-base-flow-known}
   \norm{\up(\cdot, T)}^2_{|a|} \le e^{\eta_c T} \norm{\up_0}^2_{|a|}, \quad \eta_c = \max_{x,t} \frac{\partial_t |a(x,t)|}{|a(x,t)|}. 
\end{align}
The estimate \eqref{eq:estimate-general-base-flow-known} is more optimal.
%%%
Note that if the base-flow is time-invariant, that is $\partial_t U =0 \implies \partial_t a(x,t) =0$ and $\eta_c =0$, then we have 
\begin{align}\label{eq:estimate-general-base-flow-strict}
   \norm{\up(\cdot, T)}^2_{|a|} = \norm{\up_0}^2_{|a|}. 
\end{align}
Using the equivalence of norms, the estimate \eqref{eq:estimate-general-base-flow-strict} in the weighted $L_2$-norm also holds for standard $L_2$-norm, that is 
\begin{align}
    \norm{\up(\cdot, t)}^2 \le K \norm{\up_0}^2, \quad  K = \frac{\max_{x}|a(x)|}{\min_{x}|a(x)|} > 0.
\end{align}
Therefore, if the Jacobian $a$ does not change sign within the domain then for a constant  or time-independent base-flow (\(U(x,t) = \emph{const.}\) or \(U(x,t) = U(x)\)) the $L_2$-norm of the perturbation can never grow in time. A robust numerical method should, as much as possible, replicate the energy estimates of the continuous problem. Specifically, if the continuous system does not exhibit growth, the numerical method should likewise prevent it.

\subsection{Semi-discrete analysis}
Now consider the discretization of the spatial domain $0\le x_1 < x_2 < \cdots < x_N \le L$. Let $\vb{u}:=[u_1(t), u_2(t), \cdots, u_N(t)]^T\in \mathbb{R}^N$  approximate the solution on the grid $\vb{x}:=[x_1, x_2, \cdots, x_N]^T\in \mathbb{R}^N$, that is $u_j(t) \approx u(x_j, t)$. We denote the maximum  grid-size $\Delta x = \max_{j\in \{1, 2, \cdots, N-1\}}(x_{j+1}-x_j)$ and define the discrete  $l_2([0,L])$ scalar product and norm  by
\begin{align}\label{eq:discrete-inner-product}
\langle \vb{u}, \vb{v}\rangle := \sum_{j=1}^N u_jv_j h_j, \quad h_j>0, \quad \norm{\vb{u}}^2_h := \langle \vb{u}, \vb{u}\rangle >0, \, \forall\, \vb{u} \in \mathbb{R}^N/ \{\vb{0}\},
\end{align}
where $h_j>0$ are weights of a composite quadrature rule.

We denote the semi-discrete approximation of the  nonlinear conservation law \eqref{eq:scalar-nonlinear-conservation-law} by
\begin{align}\label{eq:scalar-nonlinear-conservation-law-semi-discrete}
    \partial_t \vb{u} + \vb{P} (\vb{u}, \vb{f}, D) = 0, \quad  \vb{u}(0) = \vb{u}_0 \in \mathbb{R}^N, \quad t \in [0, T],
\end{align}
where $D$ is a consistent approximation of the first derivative on the grid, that is for a smooth $u$ we have $(D\vb{u})_j = \partial_x u|_{x = x_j} + O(\Delta x^r)$ with $r>0$. The discrete differential operator \(\vb{P} (\vb{u}, \vb{f}, D)\) denotes a consistent discrete approximation of the divergence of the PDE flux \(\partial_x f\) including the periodic boundary conditions. We introduce the notion of entropy stability for the consistent semi-discrete approximation \eqref{eq:scalar-nonlinear-conservation-law-semi-discrete}.
\begin{definition}
    A consistent semi-discrete approximation of the nonlinear conservation law \eqref{eq:scalar-nonlinear-conservation-law} is called entropy stable  (entropy conservative) if $\norm{\vb{u}(T)}^2_h \le K_n \norm{\vb{u}_0}^2_h$ ($\norm{\vb{u}(T)}^2_h = \norm{\vb{u}_0}^2_h$) for some moderate and real constant $K_n>0$ independent of the mesh parameters.
\end{definition}
%%%
High-order accurate and entropy-stable methods for nonlinear conservation laws can often be developed by combining SBP discrete derivative operators with a skew-symmetric (split-form) reformulation of the conservation laws. This approach eliminates the need to invoke the chain rule and product rule when establishing entropy stability. However, it is crucial that the linearization of the semi-discrete approximation given by \eqref{eq:scalar-nonlinear-conservation-law-semi-discrete} is strictly stable, see Definition \ref{def:strict-stability} below.

Let   \(A = \operatorname{diag}{([a(x_1,t), a(x_2, t), \cdots, a(x_N, t)]})\in \mathbb{R}^{N\times N}\) and \(Q (A, D)\in \mathbb{R}^{N\times N}\), with \(\vb{\up} \in \mathbb{R}^{N}\) denoting the approximation of the perturbation on the grid and \(Q (A, D) \vb{\up}\) denotes the linearization of the nonlinear  semi-discrete differential operator $\vb{P} (\vb{u}, \vb{f}, D)$ such that
\begin{align}\label{scalar-linear-conservation-law-discrete}
    \partial_t \vb{\up} = Q (A, D) \vb{\up}, \quad  \vb{\up}(0) = \vb{\up}_0 \in \mathbb{R}^N, \quad t \in [0, T].
\end{align}
Further, we assume that the linearized semi-discrete equation \eqref{scalar-linear-conservation-law-discrete} is a consistent approximation of the linearized conservation law \eqref{eq:scalar-linear-conservation-law}. The following definitions are significant for the current study.
\begin{definition}
    A consistent semi-discrete approximation of the linearized conservation law \eqref{eq:scalar-linear-conservation-law}  is called stable  if there exist constants $K_n, \eta_n \in \mathbb{R}$ with $K_n>0$ such that the estimate $\norm{\vb{\up}(T)}^2_h \le K_n e^{\eta_n t}\norm{\vb{\up}_0}^2_h$ holds.
\end{definition}
\noindent
If the numerical growth rate is an arbitrary real number, \(\eta_n = \eta_0 \in \mathbb{R}\), this corresponds to Lax stability \cite{LaxRichtmyer1956}. As is well known, a Lax-stable numerical method can allow non-physical numerical growth that does not occur in the continuous problem, which may require excessively fine grid resolutions to obtain accurate simulations. To be effective, the numerical growth rate \(\eta_n\) should not exceed the continuous growth rate \(\eta_c\). In particular, if the continuous problem does not support any growth, the numerical method should also exhibit no growth. This property is often referred to as strict stability; see \cite{gustafsson1995time} for further details. We now provide the formal definition.
\begin{definition}\label{def:strict-stability}
Let $\eta_c \in \mathbb{R}$ denote the continuous growth rate of the solutions of the linear conservation law \eqref{eq:scalar-linear-conservation-law}. A consistent  semi-discrete approximation of \eqref{eq:scalar-linear-conservation-law} is called strictly stable  if there exist constants $K_n, \eta_n \in \mathbb{R}$ with $K_n>0$ and $\eta_n \le \eta_c  + O(\Delta x)$,   such that the estimate $\norm{\vb{\up}(T)}^2_h \le K_ne^{\eta_n T}\norm{\vb{\up}_0}^2_h$ holds.
\end{definition}
%%%
% 
A strictly stable numerical method guarantees that the asymptotic numerical growth rate does not surpass the continuous growth rate. This is crucial for maintaining the accuracy of numerical solutions, especially in long-term simulations. For nonlinear conservation laws, strict linear stability is essential to prevent unresolved high-frequency wave modes from dominating the simulation. Consequently, the following definition is fundamental to the present study.
%%%
%%%
\begin{definition}
Consider the semi-discrete approximation \eqref{eq:scalar-nonlinear-conservation-law-semi-discrete} of the nonlinear conservation law \eqref{eq:scalar-nonlinear-conservation-law} and the linearized semi-discrete equation \eqref{scalar-linear-conservation-law-discrete}. The nonlinear semi-discrete approximation \eqref{eq:scalar-nonlinear-conservation-law-semi-discrete} is called locally energy stable if the linearized semi-discrete equation \eqref{scalar-linear-conservation-law-discrete} is strictly stable.
\end{definition}
%%%
The local energy stability properties of the nonlinear semi-discrete approximation \eqref{eq:scalar-nonlinear-conservation-law-semi-discrete} are closely related to the eigenvalues of the discrete differential operator \(Q(A,D)\) defined in \eqref{scalar-linear-conservation-law-discrete}. Specifically, let \(\lambda_{\max} := \max_{j} \Re(\lambda_j(Q))\) denote the largest real part among the eigenvalues of \(Q(A,D)\). A necessary condition for ensuring local energy stability is that \(\lambda_{\max} \leq \eta_c + O(\Delta x)\). 
Furthermore, if the discrete differential operator \(Q(A,D)\) is diagonalizable, then the estimate
\(\norm{\vb{\up}(T)}_h^2 \leq K_n e^{\lambda_{\max} T} \norm{\vb{\up}_0}_h^2\)
holds. The subsequent sections will provide more specific applications of these concepts, particularly in the context of the DP FD/DG SBP-SAT approximation of the inviscid Burgers' equation and the nonlinear SWEs.

\section{A DP SBP approximation for the scalar conservation law}\label{sec:sbp-and-scalar}
In this section, we will give a quick introduction to entropy stable DP SBP approximation for the scalar nonlinear conservation law \eqref{eq:scalar-nonlinear-conservation-law}. Detailed analysis of local linear/energy stability for the approximation will be presented, focusing on the Burgers' equation. For systems of conservation laws and more elaborate discuss, we refer the reader to the recent work \cite{duru2024dualpairingsummationbypartsfinitedifference}.
\subsection{Skew-symmetric split-form}\label{sec:split-form}
Entropy stable schemes for nonlinear conservation laws typically begins at the continuous level by deriving the skew-symmetric form.
The skew-symmetric form is often found by splitting the divergence of the flux as a convex combination of the nonlinear conservation laws in flux/conservative form and the quasi-linear form  so that the energy method can be applied without the chain rule and/or the product rule. For a scalar PDE flux $f(u)$ we have
$$
\partial_x f(u)\equiv F(u, f, \partial_x):=  \alpha \partial_x f +  (1-\alpha)  \partial_u f \partial_x u, \quad \alpha \in [0, 1],
$$
 which is a convex combination of the flux/conservative form and the advective form. For Burgers' equation with $f(u) = u^2/2$ we have $\alpha = 2/3$ and
 $
F(u, f, \partial_x):=    1/3\partial_x (u^2) +  1/3  u \partial_x u.
$
It is significantly noteworthy that if $\alpha = 1$ we regain the flux/conservative form, that is $F(u, f, \partial_x) =    \partial_x f$.
%%%%
\subsection{The DP SBP  operator for {$d/dx$}}\label{sec:SBP framework}
Recently, the  DP SBP FD framework \cite{MATTSSON2017upwindsbp,DovgilovichSofronov2015,williams2024drp} was introduced to improve the accuracy and flexibility of numerical approximations of PDEs. The  DP SBP  operators are a pair of forward and backward finite difference stencils that together obey  the SBP property. The DP upwind SBP framework for the first derivative $d/dx$ \cite{MATTSSON2017upwindsbp,DovgilovichSofronov2015} can be expressed through the following assumptions:
\begin{enumerate}
    \item[A.1] There exists $H: \mathbb{R}^N \to \mathbb{R}^N$ which defines a positive discrete measure
    \begin{align*}
        \langle\mathbf{g}, \mathbf{g} \rangle_{H}
            = \mathbf{g}^T{H}\mathbf{g} > 0,
            \quad \forall \, \mathbf{g} \in \mathbb{R}^N, 
        \quad\langle \mathbf{1}, \mathbf{g} \rangle_{H}
            = \sum_{j = 1}^N h_j  g(x_j) \to \int_{0}^{L} g(x) \dd{x},
            %\quad \forall\, g \in L^2(\Omega).
    \end{align*}
     for all $g \in C^0([0, L])$.
%     %
    \item [A.2] There exists a pair of linear operators $D_{\pm}: \mathbb{R}^N \to \mathbb{R}^N $ with
    $\left(D_{\pm} \mathbf{f}\right)_j = \dv{f}{x}|_{x=x_j}$,
    for all $j \in \{1, 2, \dots, N\}$  and $f \in V^p$ where $V^p$ is a polynomial space of at most degree $p \ge 0$.
    \item [A.3] The linear operators $D_{\pm}$ together obey $\langle D_{+} \mathbf{f}, \mathbf{g} \rangle_{H} + \langle \mathbf{f}, D_{-} \mathbf{g} \rangle_{H} =f_Ng_N - f_1g_1$ for all $\mathbf{f}, \, \mathbf{g} \in \mathbb{R}^N$.

    \item [A.4] The linear operators $D_{\pm}$ together obey 
    %\begin{align}\label{eq:B4}
         $\langle \mathbf{f}, \left(D_{+}-D_{-} \right)\mathbf{f} \rangle_{H} \le 0$ for all $\mathbf{f} \in \mathbb{R}^N$.
    %\end{align}
    %for all $\mathbf{f} \in \mathbb{R}^n$.
    %
\end{enumerate}
%
%%%
\begin{remark}
    The DP SBP DG framework of \cite{GLAUBITZ2025113841} satisfies Assumptions ($\mathbf{A}$.1)--($\mathbf{A}$.4) with $x_j$, $w_j>0$ being the nodes and weights of a Legendre-Gauss-Lobatto (LGL) numerical quadrature.
\end{remark}
\begin{remark}
The traditional SBP operator $(D, H)$--based on central FD~\cite{kreiss1974finite,BStrand1994,fernandez2014review,svard2014review}--or collocated DGSEM~\cite{gassner2013skew,taylor2010compatible,fluxo_SplitForm}--also satisfies the DP-SBP framework, that is Assumptions ($\mathbf{A}$.1)--($\mathbf{A}$.4), with   $D_- =D_+ =D$  and $\langle \mathbf{f}, \left(D_{+}-D_{-} \right)\mathbf{f} \rangle_{H} = 0$.
 Similarly, given the upwind DP-SBP operator $(D_-, D_+, H)$ which satisfies Assumptions ($\mathbf{A}$.1)--($\mathbf{A}$.4), the averaged operator
$
D := \frac12 \left(D_+ + D_-\right)
$
is a discrete derivative operator which satisfies the traditional SBP framework $(D, H)$. 
\end{remark}

The SBP operators alone do not guarantee numerical stability; they must be applied carefully with appropriate boundary treatments. The DP-SBP framework offers advantages over traditional SBP schemes by providing pairs of derivative operators and increased flexibility, enabling the design of schemes with desirable properties. For conservative PDEs, energy/entropy conservation can be achieved due to the DP-SBP property (Assumption $\mathbf{A}.3$), where operators $\mathcal{Q}_{\pm} :=  H{D}_{\pm}  - \frac{1}{2}B$, with $B = \diag([-1, 0, \cdots, 0, 1])$, satisfy $\mathcal{Q}_{-} +\mathcal{Q}_{+}^T =0.$ Additionally, entropy or energy dissipation schemes can be constructed, useful for upwinding or damping unresolved modes, primarily enabled by the upwind DP-SBP (Assumption $\mathbf{A}.4$). 
Using penalties we will  extend the DP SBP operators to multi-block operators and implement the periodic boundary conditions weakly.

%\subsection{Multi-block SBP operators}
For simplicity, we split the spatial domain into two elements, with $K=2$ and
$\Omega = \Omega_k \cup \Omega_{k+1}$. The method and analysis extend to many more elements. For each element $\Omega_k$, $k \in \{1, 2\}$, we discretize the domain into $N$ grid points and sample the solutions on the grids. For the multi-block FD method the grids are equidistant nodes and for the DG method the grids are given by the LGL nodes. The semi-discrete solution is denoted $\mathbf{u}^k =  ({u}_1^k, {u}_2^k, \cdots, {u}_N^k)^T$. The two elements are connected at the interface $x = x_k$ with the interface conditions
\begin{align*}
   \lJump{\vb{u}} \rJump:  = {u}_1^{k+1}-{u}_N^k= 0, \quad \lJump{\vb{f}} \rJump:  = {f}_1^{k+1}-{f}_N^k= 0. 
\end{align*}
%%%
For the global solutions in the two elements, we introduce the  augmented  discrete solution vector and discrete flux vector, given by 
\begin{align*}
\mathbf{u} = \begin{bmatrix}
{u}_1^k, \cdots, {u}_N^k,{u}_1^{k+1},  \cdots, {u}_N^{k+1}
\end{bmatrix}^T \in \mathbb{R}^{2N}, \quad 
\mathbf{f} = \begin{bmatrix}
{f}_1^k, \cdots, {f}_N^k,{f}_1^{k+1},  \cdots, \mathbf{f}_N^{k+1}
\end{bmatrix}^T  \in \mathbb{R}^{2N}.    
\end{align*}
%%%
To implement the interface conditions and the periodic boundary conditions we also define the boundary and interface penalty matrices
{%\small
\scriptsize
\begin{align}\label{eq:SBP_SAT_and_Interface_Matrices}
\vec{B}_{I} = \begin{pmatrix}
-\vb{e}_N\vb{e}_N^T & \vb{e}_N\vb{e}_1^T\\
-\vb{e}_1\vb{e}_N^T & \vb{e}_1\vb{e}_1^T
\end{pmatrix},
\quad
\widetilde{\vec{B}}_{I} = \begin{pmatrix}
-\vb{e}_N\vb{e}_N^T & \vb{e}_N\vb{e}_1^T\\
\vb{e}_1\vb{e}_n^T & -\vb{e}_1\vb{e}_1^T
\end{pmatrix}, \quad \vec{B}_{N} = \begin{pmatrix}
\vb{e}_1\vb{e}_1^T & -\vb{e}_1\vb{e}_N^T\\
\vb{e}_N\vb{e}_1^T & - \vb{e}_N\vb{e}_N^T
\end{pmatrix},
\quad
\widetilde{\vec{B}}_{N} = \begin{pmatrix}
-\vb{e}_1\vb{e}_1^T & \vb{e}_1\vb{e}_N^T\\
\vb{e}_N\vb{e}_1^T & - \vb{e}_N\vb{e}_N^T
\end{pmatrix},
\end{align}
}
where $ \vb{e}_1^T \vb{u} = u(x_1)$ and $ \vb{e}_N^T \vb{u} = u(x_N)$.
We have
{
$$
\vec{B}_{I}\mathbf{u} = \left(
0, \cdots, 0,\lJump{\vb{u}}\rJump, \lJump{\vb{u}}\rJump, 0,  \cdots, 0
\right)^T, \quad \widetilde{\vec{B}}_{I}\mathbf{u} = \left(
0, \cdots, 0,\lJump{\vb{u}}\rJump, -\lJump{\vb{u}}\rJump, 0,  \cdots, 0
\right)^T, 
$$
The interface matrices $\vec{B}_{I}$, $\widetilde{\vec{B}}_{I}$ will be used to weakly couple the elements together to the global domain and the boundary matrices $\vec{B}_{N}$, $\widetilde{\vec{B}}_{N}$ will weakly implement the period boundary conditions. Note also that the matrices $\widetilde{\vec{B}}_{I}$, $\widetilde{\vec{B}}_{N}$ are symmetric and negative semi-definite, that is $\widetilde{\vec{B}}_{I} = \widetilde{\vec{B}}_{I}^T$, $\widetilde{\vec{B}}_{N} = \widetilde{\vec{B}}_{N}^T$ and 
$
\mathbf{u}^T\widetilde{\vec{B}}_{I}\mathbf{u} = -\lJump{\vb{u}}\rJump^2=-(u_1^{k+1}-u_N^{k})^2\le 0,
$
$
\mathbf{u}^T\widetilde{\vec{B}}_{N}\mathbf{u} = -(u_1^{1}-u_N^{K})^2\le 0.
$
We introduce the multi-block operators
{\small
\begin{align}\label{eq:periodic_interface_SBP_SAT}
\vec{H} = \begin{pmatrix}
H & 0\\
0 & H
\end{pmatrix},
\quad
\vec{D} = \begin{pmatrix}
D & 0\\
0 & D
\end{pmatrix},
\quad
\vec{D}_{\pm} = \begin{pmatrix}
D_{\pm} & 0\\
0 & D_{\pm}
\end{pmatrix},
\quad 
\widetilde{\vec{D}}_{\pm} = \vec{D}_{\pm} + \frac12 \vec{H}^{-1}\left(\vec{B}_N + \vec{B}_I\right),
\end{align}
%%%
 where the penalty terms  weakly implement the interface and  the periodic boundary conditions. The penalized multi-block operators $\widetilde{\vec{D}}_{\pm}$ also satisfy the identities
 %\vspace{-0.25cm}
 {
\begin{align}\label{eq:B4_periodic}
\inp{\widetilde{\vec{D}}_{+} \mathbf{f}}{\mathbf{g}}_{\vec{H}} + \inp{ \mathbf{f}}{\widetilde{\vec{D}}_{-} \mathbf{g}}_{\vec{H}} = 0,
    \,
\inp{ \mathbf{f}}{\left(\widetilde{\vec{D}}_{+}-\widetilde{\vec{D}}_{-} \right)\mathbf{f}}_{\vec{H}}=    \inp{\mathbf{f}} {\left(\vec{D}_{+}-\vec{D}_{-} \right)\mathbf{f}}_{\vec{H}} \le 0, \quad \forall \, \mathbf{f}, \, \mathbf{g} \in \mathbb{R}^{2N}.
\end{align}
}
%%%
 For simplicity, we will ignore the $\sim$ in the penalized SBP operators and write $\vec{D}_{\pm}=\widetilde{\vec{D}}_{\pm}$ where $\widetilde{\vec{D}}_{\pm}$ are multi-block and  periodic operators implemented using penalties as in \eqref{eq:periodic_interface_SBP_SAT}.

%%%%
%%%%
\subsection{An entropy stable DP SBP approximation}
A skew-symmetric SBP approximation of the divergence of the PDE flux including the periodic boundary condition is
 \begin{equation}
	\vb{F} (\vb{u}, \vb{f}, \vec{D}):= \alpha \vb{D} \vb{f} + (1-\alpha) \vb{\partial_u f} \vb{D} \vb{u}, \quad  \vec{D} = \frac{\vec{D}_{+}+\vec{D}_{-}}{2}.
\end{equation}
For Burgers' equation we have  
 $
	\vb{F} (\vb{u}, \vb{f}, \vec{D}) := 1/3 \vb{D} \left(\vb{u^2}\right) + 1/3 \vb{ u} \vb{D} \vb{u}.
$
 A multi-block entropy stable DP SBP semi-discrete   approximation of the systems of 1D conservation laws \eqref{eq:scalar-nonlinear-conservation-law} with volume upwinding reads
 \begin{align}\label{eq:skew_symm_upwind_SBP_SAT-multi-block}
    \dv{}{t} \mathbf{u} + \vb{P} (\vb{u}, \vb{f}, \vec{D}) = 0,\quad \vb{P} (\vb{u}, \vb{f}, \vec{D}):=\vb{F}(\vb u, \vb{f}(\vb{u}), \vec{D})-\gamma \vb{\Di} \vb{g}, \quad \vb{\Di} = \frac{\vec{D}_{+}-\vec{D}_{-}}{2}, 
\end{align}
where $\gamma >0$ is the volume upwind parameter and $\vb{g}$ is the entropy variable. If the continuous problem \eqref{eq:scalar-nonlinear-conservation-law} admits a square entropy and \(\vec{D}_{+},\vec{D}_{-}\) satisfy the DP SBP properties then it follows that the nonlinear discrete energy/entropy estimate holds, that is 
\begin{equation}
    \dv{}{t}{\norm{\vb{u}(t)}^2_h} = \gamma \vb{g}^\top (\vb{H}\vb{\Di} )\vb{g} \le 0 \implies \norm{\vb{u}(T)}^2_h \le \norm{\vb{u}(0)}^2_h.
\end{equation}
In the absence of volume DP upwinding, that is $\gamma =0$ or $ \vb{\Di}=0$, the total energy/entropy is conserved, $\norm{\vb{u}(T)}^2_h = \norm{\vb{u}(0)}^2_h$. An entropy/nonlinearly stable numerical method can ensure robustness and the boundedness  of the numerical solution for nonlinear problem. We refer the reader to \cite{duru2024dualpairingsummationbypartsfinitedifference} for more elaborate discussions. However, as shown in \cite{gassner2022stability}, a numerical method can be nonlinearly stable and at the same allow the growth of unresolved nonlinear numerical wave modes which can corrupt the results of numerical simulations of nonlinear conservation laws. 
%%
%%%
\subsection{Local linear/energy stability}
%%%
Our primary goal in this study is the investigation of the linear stability properties of the recent paper \cite{duru2024dualpairingsummationbypartsfinitedifference}. As a model problem we will consider the model problem \eqref{eq:scalar-nonlinear-conservation-law} and the corresponding DP-SBP entropy stable approximation \eqref{eq:skew_symm_upwind_SBP_SAT-multi-block}.

Next, consider a smooth base-flow \(\vb{U}\) and  linearize the nonlinear semi-discrete approximation \eqref{eq:skew_symm_upwind_SBP_SAT-multi-block} 
\begin{equation}\label{eq:linear_upwind_SBP_SAT-multi-block} 
    \dv{}{t}\vb{\up} = \vb{Q}\vb{\up}, \quad \vb{Q}\vb{\up}=
    -  \vb{D} (\vb{A} \vb{\up})
    - \underbrace{(1 - \alpha) \left(\left( \vb{B}  \vb{D}_U \vb{\up} + \vb{A} \vb{D} \vb{\up}\right) - \vb{D} \left(\vb{A} \vb{\up}\right)\right)}_{\text{destablizing term}}
    + \gamma \boldsymbol{\Di} ( \vb{C} \vb{\up}).
\end{equation}
Here $\vb{D}_U = \operatorname{diag}{([(\vb{D}\vb{U})_1, (\vb{D}\vb{U})_2, \cdots (\vb{D}\vb{U})_{2N}])}$, and $A,B,C$
are diagonal matrices, with $A_{ii}=a(x_i,t)$, $B_{ii}=b(x_i,t)$, $C_{ii}=c(x_i,t)$, where
\( a = \partial_uf|_{u =U}\) is the Jacobian of the flux, \( b = \partial_{U}a\) is the Hessian of the flux, and  \( c = \partial_ug|_{u =U}\) is the Jacobian of the entropy variable with respect to the conserved variable. For Burgers' equation \( a =  U\), \( b =  1\), and \( c =  1\), yielding
\begin{equation}\label{eq:linear_upwind_SBP_SAT-multi-block-Burgers} 
    \dv{}{t}\vb{\up} = \vb{Q}\vb{\up}, \quad 
    \vb{Q} =  -\vb{D} \vb{A} 
    - \underbrace{(1 - \alpha) \left( \vb{D}_U + \vb{A} \vb{D} -\vb{D} \vb{A}\right)}_{\text{destabilizing term}} +
     \gamma \boldsymbol{\Di}, 
\end{equation}
In \eqref{eq:linear_upwind_SBP_SAT-multi-block} and  \eqref{eq:linear_upwind_SBP_SAT-multi-block}, the "destabilizing" terms inside the bracket with the factor\((1-\alpha)\) are a residual of  the chain/product rule. That is 
\begin{align*}
   \vb{B}  \vb{D}_{U} \vb{\up} + \vb{A} \vb{D} \vb{\up} -\vb{D} (\vb{A} \vb{\up}) &\approx (\partial_Ua) (\partial_x U)  \up + a (\partial_x \up) - \partial_x(a \up)\\
   &=  (\partial_x a)  \up + a (\partial_x \up) - \partial_x(a \up) \equiv 0. 
\end{align*}
%%%
Note that for Burgers' equation $a=U$.
The "destabilizing" terms will vanish if \(\alpha = 1\) or if the discrete derivative operator \(\mathbf{D}\) satisfies both the chain rule and the product rule.
Note in particular that \(\alpha = 1\) and $\gamma >0$ will correspond  to the upwind method of \cite{GLAUBITZ2025113841,ranocha2023highorder}, which are well-known to be linearly stable. However, these schemes are not nonlinearly stable and may lack the robustness necessary for simulating nonlinear phenomena such as turbulence and shocks, see \cite{duru2024dualpairingsummationbypartsfinitedifference}.
When $\alpha \ne 1$ and $\gamma =0$, to ensure local linear stability the discrete derivative operator \(\mathbf{D}\) must satisfy both the chain rule and the product rule. However, \(\mathbf{D}\) will only satisfy these rules in the trivial case of a constant background flow \(\mathbf{U} = \text{const.}\) It appears that the skew-symmetric flux splitting, which ensures nonlinear and entropy stability, can actually destabilize the linearized semi-discrete approximation and amplify grid-scale errors. In the context of finite element discretization \cite{Spiegel2015,Manzanero2018}, efforts such as over-integration have been employed to reduce aliasing errors. Nonetheless, these approaches often do not yield satisfactory results, as they can introduce other stability issues.

In the following analysis, 
as in the continuous setting, we will consider the background flow such that the coefficient $a(x,t)$ does not change sign within the domain, and the linearize continuous problem does not support growth. Thus if $\alpha = 1$ and $\gamma =0$, we have
\begin{align}\label{eq:estimate-general-base-flow-known-discrete}
   \norm{\vb{\up}(\cdot, T)}^2_{|AH|} \le e^{\eta_d T} \norm{\vb{\up}_0}^2_{|AH|}, \quad \eta_d = \max_{j,t} \frac{\partial_t |a(x_j,t)|}{|a(x_j,t)|}. 
\end{align}
%%%
%%%
Therefore, if  $\alpha = 1$ and $\gamma =0$ and  in addition  the base-flow is time-invariant that is $\partial_t U =0 \implies \partial_t a(x,t) =0 \implies\eta_d =0$, then we have 
%%%%%
\begin{align}\label{eq:estimate-general-base-flow-strict-discrete}
    \norm{\vb{\up}(\cdot, T)}^2_{H} \le K \norm{\vb{\up}_0}^2_{H},  \quad  K = \frac{\max_{j}|a(x_j)|}{\min_{j}|a(x_j)|} > 0. 
\end{align}
%%%%
Note that in \eqref{eq:estimate-general-base-flow-known-discrete} and \eqref{eq:estimate-general-base-flow-strict-discrete}, the numerical growth rate is bounded by the continuous growth rate, i.e., \(\eta_d \leq \eta_c\), thus yield local linear energy stability fot the discretization. However, setting \(\alpha = 1\) does not yield an entropy-stable scheme for the nonlinear problem. Furthermore, detailed numerical experiments conducted in \cite{duru2024dualpairingsummationbypartsfinitedifference} demonstrate that the linearly stable scheme with \(\alpha = 1\) is not robust; when nonlinear phenomena are dominant it crashes before reaching the final time, even with a non-vanishing upwind parameter (\(\gamma > 0\)). Our goal is to show that the DP FD/DG SBP method \eqref{eq:skew_symm_upwind_SBP_SAT-multi-block} can  simultaneously preserve both nonlinear/entropy stability and local linear stability. Achieving this ensures robustness and prevents grid-scale errors from contaminating the numerical simulations.

For the Burgers' equation, discretized using the DP FD/DG SBP method \eqref{eq:skew_symm_upwind_SBP_SAT-multi-block}, we will demonstrate that, for \(\alpha = {2}/{3}\), it is possible to select a positive volume upwinding parameter \(\gamma > 0\) such that the scheme is both entropy stable and locally linear energy stable. This result will ensure both entropy/nonlinear stability and local linear/energy stability, thereby improving the accuracy of numerical simulations and enabling the proof of convergence for sufficiently smooth solutions.
As discussed in the previous section, the local linear and energy stability of the semi-discrete approximation \eqref{eq:skew_symm_upwind_SBP_SAT-multi-block} for the nonlinear conservation law is determined by the eigenvalues of the linear operator \(\vb{Q}\), defined in \eqref{eq:linear_upwind_SBP_SAT-multi-block-Burgers}. In particular, if the base-flow is sufficiently smooth and does not change sign within the domain, then we must have \(\Re{\lambda(\vb{Q})} \le 0\).
\subsection{A finite volume reformulation  for the Burgers' equation}
Here, we will consider specifically the Burgers' equation and the semi-discrete DP SBP method \eqref{eq:linear_upwind_SBP_SAT-multi-block} in a single element and reformulate the reformulate as finite volume method by introducing numerical fluxes $\vb{f}_{i\pm\frac12}^N$. For the Burgers' equation, with the PDE flux $f(u) =u^2/2$, the finite volume reformulation of the DP FD SBP semi-discrete approximation \eqref{eq:skew_symm_upwind_SBP_SAT-multi-block} is
\begin{align}\label{eq:skew_symm_upwind_SBP_SAT-finite-volume}
    \dv{}{t}{\vb{u}_i} + \frac{\vb{f}_{i+\frac12}^N - \vb{f}_{i-\frac12}^N}{\Delta{x}}  = 0.
\end{align}
The numerical fluxes $\vb{f}_{i\pm\frac12}^N$ decompose into 
%\begin{align}
$
    \vb{f}_{i\pm\frac12}^N= \vb{f}_{i\pm\frac12}^{C} + \vb{f}_{i\pm\frac12}^{R} + \vb{f}_{i\pm\frac12}^{U},
    $
%\end{align}
where $\vb{f}_{i\pm \frac12}^{C}$ is the centered flux, $\vb{f}_{i\pm \frac12}^{R}$ is the residual flux, $\vb{f}_{i\pm \frac12}^{U}$ is the upwind flux, and are defined by
\begin{align*}
	\vb{f}_{i+\frac12}^{C}
		&= \Delta{x}\sum_{j > i} D_{ij} \vb{f}_j + \frac12 \vb{f}_i, \quad
	\vb{f}_{i-\frac12}^{C}
		= -\Delta{x}\sum_{j < i} D_{ij} \vb{f}_j + \frac12 \vb{f}_i, \\
        %\quad \beta = 1 - \Delta{x}\sum_{j>i} D_{ij}=1 + \Delta{x}\sum_{j < i} D_{ij}, \\
	\vb{f}_{i+\frac12}^{R}
        &= -\frac{1-\alpha}{2} \Delta{x}\sum_{j>i} D_{ij} (\vb{u}_i - \vb{u}_j)^2, \quad
	\vb{f}_{i-\frac12}^{R}
        = \frac{1-\alpha}{2} \Delta{x}\sum_{j<i} D_{ij} (\vb{u}_i - \vb{u}_j)^2, \\
	\vb{f}_{i+\frac12}^{U}
		&= -\gamma\Delta{x}\left(\frac{1}{2} \Di_{ii} \vb{u}_i + \sum_{j>i} \Di_{ij} \vb{u}_j\right), \quad
	\vb{f}_{i-\frac12}^{U}
		= \gamma \Delta{x}\left(\frac{1}{2} \Di_{ii} \vb{u}_i + \sum_{j<i} \Di_{ij} \vb{u}_j\right).
\end{align*}
Note that when $\alpha = 1$ and $\gamma =0$ we get \(\vb{f}_{i\pm\frac12}^{R} =0, \vb{f}_{i\pm\frac12}^{U}=0\) which yields the  SBP scheme applied to the flux/conservative form.  With $\alpha = 2/3$ and $\gamma = 0$ or $\vb{\Di} =0$ we obtain \( \vb{f}_{i\pm\frac12}^{U}=0\) which corresponds to the skew-symmetric and entropy conserving SBP scheme without volume upwinding. However, $\alpha = 2/3$ and $\gamma > 0$ yields the skew-symmetric and entropy stable DP-SBP scheme with volume upwinding.

Next we linearized the semi-discrete approximation \eqref{eq:skew_symm_upwind_SBP_SAT-finite-volume} of the Burgers' equation around a smooth base-flow, $\vb{U}_{i} = U(x_i, t)$, we have
\begin{align}\label{eq:skew_symm_upwind_SBP_SAT-finite-volume-linear}
    \dv{}{t}{\delta \vb{u}_i} + \frac{\delta\vb{f}_{i+\frac12}^N - \delta\vb{f}_{i-\frac12}^N}{\Delta{x}}  = 0.
\end{align}
Again, the linear numerical  fluxes $\delta\vb{f}_{i\pm\frac12}^N$ decompose into 
\begin{align}
    \delta\vb{f}_{i\pm\frac12}^N= \delta\vb{f}_{i\pm\frac12}^{C} + \delta\vb{f}_{i\pm\frac12}^{R} + \delta\vb{f}_{i\pm\frac12}^{U},
\end{align}
where $\delta\vb{f}_{i\pm \frac12}^{C}$ is the centered flux, $\delta\vb{f}_{i\pm \frac12}^{R}$ is the residual flux, $\delta\vb{f}_{i\pm \frac12}^{U}$ is the upwind flux, and are defined by
{%\footnotesize
\scriptsize
\begin{align*}
	\delta\vb{f}_{i+\frac12}^{C}
		&= \Delta{x}\sum_{j > i} D_{ij} \vb{\up}_j\vb{U}_j +  \vb{\up}_i\vb{U}_i, \quad
	\delta\vb{f}_{i-\frac12}^{C}
		= -\Delta{x}\sum_{j < i} D_{ij} \vb{\up}_j\vb{U}_j +  \vb{\up}_i\vb{U}_i, %\quad \beta = 1 - \Delta{x}\sum_{j>i} D_{ij}=1 + \Delta{x}\sum_{j < i} D_{ij}, 
        \\
	\delta\vb{f}_{i+\frac12}^{R}
        &= -\left({1-\alpha}\right) \Delta{x}\sum_{j>i} D_{ij} (\delta \vb{u}_i - \delta\vb{u}_j)( \vb{U}_i - \vb{U}_j), \quad
	\delta\vb{f}_{i-\frac12}^{R}
        = \left({1-\alpha}\right) \Delta{x}\sum_{j<i} D_{ij} (\delta \vb{u}_i - \delta\vb{u}_j)( \vb{U}_i - \vb{U}_j), \\
	\delta\vb{f}_{i+\frac12}^{U}
		&= -\gamma\Delta{x}\left(\frac{1}{2} \Di_{ii} \delta\vb{u}_i + \sum_{j>i} \Di_{ij} \delta\vb{u}_j\right), \quad
	\delta\vb{f}_{i-\frac12}^{U}
		= \gamma \Delta{x}\left(\frac{1}{2} \Di_{ii} \delta\vb{u}_i + \sum_{j<i} \Di_{ij} \delta\vb{u}_j\right).
\end{align*}
}
Recall that the linear centered flux \(\delta \mathbf{f}_{i+\frac{1}{2}}^{C}\) is energy conservative. In contrast, the linear residual flux \(\delta \mathbf{f}_{i+\frac{1}{2}}^{R}\) is not sign definite and can potentially lead to energy growth when \(\alpha \neq 1\). By construction, the linear upwind  flux \(\delta \mathbf{f}_{i+\frac{1}{2}}^{U}\) is energy dissipative. 
Our goal is to select the upwind parameter \(\gamma > 0\) -- sufficiently large and optimally chosen -- so that the linear residual flux \(\delta \mathbf{f}_{i+\frac{1}{2}}^{R}\) is effectively "absorbed" by the linear upwind  flux. In other words, we want:
\(\delta \mathbf{f}_{i+\frac{1}{2}}^{N} = \delta \mathbf{f}_{i+\frac{1}{2}}^{C} + \bar{\gamma} \delta \mathbf{f}_{i+\frac{1}{2}}^{U},\)
for some \(\bar{\gamma} > 0\). 
The choice of \(\gamma > 0\) or \(\bar{\gamma} > 0\) depends on the order of accuracy of the specific DP-SBP operator employed. For first-order accurate DP-SBP operators \(D_{\pm}\), the interior stencils are
 \begin{align}
            \Delta{x}\Diff &= \frac 1 {2} \left[-1, 0, 1\right], \quad \Delta{x} \Di = \frac 1 {2} \left[1, -2, 1\right], 
\end{align}
and the linear flux
\begin{align}
    2\vb{\fp}_{i+\frac 1 2}^N
                &= \underbrace{\left(\up_{i+1} U_{i+1} + \up_i U_i\right)}_{2\vb{\fp}_{i+\frac 1 2}^C}  
                + \bar{\gamma} \underbrace{\left(\up_{i} - \up_{i+1}\right)}_{2\vb{\fp}_{i+\frac 1 2}^U},
\end{align}
where $\bar{\gamma} = (1-\alpha) (U_{i+1} - U_i) + \gamma$. Thus if 
$ \gamma \ge \gamma_{opt}=(1-\alpha) \max_{i}|(U_{i+1} - U_i)| \approx (1-\alpha)\max_{i}|\Delta x \partial_xU(x_i,t)|$ then we must have $\bar{\gamma}>0$. Note that for $\alpha = 2/3$, we have $\gamma_{\mathrm{opt}} = 1/3\max_{i}|\Delta x \partial_xU(x_i,t)|  \le  2/3\max_{i}|U_{i}| < \gamma_{\mathrm{LF}}= \max_{i}| U_{i}|$ where $\gamma_{\mathrm{LF}}$ denotes the upwind parameter given by the global Lax-Friederichs' scheme \cite{Lax-Friedrichs1971}. For a sufficiently smooth base-flow, we have $\gamma_{\mathrm{opt}} =  1/3\max_{i}|\Delta x \partial_xU(x_i,t)| \to 0$ as $\Delta{x} \to 0$. 

For the third order accurate DP-SBP operators $D_{\pm}$ we have the interior stencils
 \begin{align}
            \Delta{x}\Diff &= \frac 1 {12} \left[1, -8, 0, 8, -1\right] , \quad \Delta{x} \Di = \frac 1 {12} \left[-1, 4, -6, 4, -1\right], 
\end{align}
and the linear flux
\begin{align*}
12\vb{\fp}_{i+\frac 1 2}
                &= \left(
                    -\up_{i+2} U_{i+2}
                    + 8\up_{i+1} U_{i+1}
                    + 5\up_i U_i\right) \\
                    &+ (1-\alpha) \left( 8(U_{i+1} - U_i) (\up_{i} - \up_{i+1})
                    -(U_{i+2} - U_i) (\up_{i} - \up_{i+2})
                \right)\\
                &+ \gamma \left(4\left(\up_{i} - \up_{i+1}\right) -\left(\up_{i} -\up_{i+2}\right)\right) .
\end{align*}
Thus for a smooth base-flow $U$, with $U_{i+1} - U_i \approx \Delta x \partial_x U(x_i,t)$ and $U_{i+2} - U_i \approx 2\Delta x \partial_x U(x_i,t)$, we have
\begin{align*}
12\vb{\fp}_{i+\frac 1 2}
                &= \underbrace{\left(
                    -\up_{i+2} U_{i+2}
                    + 8\up_{i+1} U_{i+1}
                    + 5\up_i U_i\right)}_{12\vb{\fp}_{i+\frac 1 2}^C} 
                    + \bar{\gamma} \underbrace{\left(4\left(\up_{i} - \up_{i+1}\right) -\left(\up_{i} -\up_{i+2}\right)\right)}_{12\vb{\fp}_{i+\frac 1 2}^U}, 
\end{align*}
where $\bar{\gamma} = 2(1-\alpha)\Delta x \partial_x U(x_i,t) + \gamma$. Thus if 
$ \gamma \ge \gamma_{opt}=(1-\alpha) \max_{i}|2\Delta x \partial_x U(x_i,t)|$ then we must have $\bar{\gamma}\ge 0$. Note again that for a sufficiently smooth base-flow, we have $\gamma_{\mathrm{opt}} =  (1-\alpha) \max_{i}|2\Delta x \partial_x U(x_i)| \to 0$ as $\Delta{x} \to 0$. 
We have extended the analysis to higher-order accurate DP FD SBP, for odd order (1, 3, 5, 7, 9) accurate operators, and have extrapolated the analysis to even order (2, 4, 6, 8, 10) accurate operators. The results are summarized in Table \ref{tab:Opt-upwind-param-Burgers} below. Numerical experiments performed in the next section show that the results also hold for the DP DRP SBP operators \cite{williams2024drp} and the DP DG SBP operators \cite{GLAUBITZ2025113841}.
\begin{remark}
    We remark that in deriving the optimal upwind parameters $\gamma_{\mathrm{opt}}$ given in in Table \ref{tab:Opt-upwind-param-Burgers}, we have made necessary smoothness assumption on the base-flow $\vb{U}$ only. No smoothness assumption is given for the perturbation $\vb{\up}$ which can be an unresolved numerical mode that lives on the grid. 
\end{remark}
\begin{table}[htbp]
    \centering
    \begin{subcaptionblock}{0.49\textwidth}
        \centering
        \begin{tabular}{lcl}
            \toprule
            order && bound for interior local linear stability \\
            \midrule
            1 && \(\gamma_{\mathrm{opt}} = (1-\alpha) \max_{i}\abs{\Delta x \, \partial_x U(x_i)}\)  \\
            3 && \(\gamma_{\mathrm{opt}} = (1-\alpha) \max_{i}\abs{2\Delta x \, \partial_x U(x_i)}\)  \\
            5 && \(\gamma_{\mathrm{opt}} = (1-\alpha) \max_{i}\abs{3\Delta x \, \partial_x U(x_i)}\)  \\
            7 && \(\gamma_{\mathrm{opt}} = (1-\alpha) \max_{i}\abs{4\Delta x \, \partial_x U(x_i)}\)  \\
            9 && \(\gamma_{\mathrm{opt}} = (1-\alpha) \max_{i}\abs{5\Delta x \, \partial_x U(x_i)}\)  \\
            \bottomrule
        \end{tabular}
        \caption{Odd interior order accurate operators}
    \end{subcaptionblock}
    \begin{subcaptionblock}{0.49\textwidth}
        \centering
        \begin{tabular}{lcl}
            \toprule
            order && bound for interior local linear stability \\
            \midrule
            2 && \(\gamma_{\mathrm{opt}} = (1-\alpha) \max_{i}\abs{\frac32\Delta x \, \partial_x U(x_i)}\)  \\
            4 && \(\gamma_{\mathrm{opt}} = (1-\alpha) \max_{i}\abs{\frac52\Delta x \, \partial_x U(x_i)}\)  \\
            6 && \(\gamma_{\mathrm{opt}} = (1-\alpha) \max_{i}\abs{\frac72\Delta x \, \partial_x U(x_i)}\)  \\
            8 && \(\gamma_{\mathrm{opt}} = (1-\alpha) \max_{i}\abs{\frac92\Delta x \, \partial_x U(x_i)}\)  \\
            10 &&\(\gamma_{\mathrm{opt}} = (1-\alpha) \max_{i}\abs{\frac{11}{2}\Delta x \, \partial_x U(x_i)}\)  \\
            \bottomrule
        \end{tabular}
        \caption{Even interior order accurate operators}
    \end{subcaptionblock}
    \caption{Optimal upwind parameters for the Burgers' equation }
    \label{tab:Opt-upwind-param-Burgers}
\end{table}

\section{Some instructive numerical experiments: Burgers' equation}\label{sec:instructive-numerical-experiments}
We will now conduct some instructive numerical experiments that illustrate and fully align with the analyses presented in the previous section. Specifically, we consider the Burgers equation on the unit interval \(x \in [0, 1]\), with the initial condition
\(u(x, 0) = 2 + \sin(\pi (x - 0.7)),\)
and periodic boundary conditions, \(u(0, t) = u(1, t)\) for all \(t \geq 0\). Note that the initial condition is strictly positive for all  \(x \in [0, 1]\), i.e., \(u(x, 0) > 0\), and does not change sign within the domain  \(x \in [0, 1]\). Consequently, at the continuous level, the linearization of Burgers' equation around this initial condition, \(U(x) = u(x, 0)\), does not support the growth of perturbations. A good numerical method, when linearized around the base-flow, should not support growth of perturbations.

We have linearized the the DP SBP approximation \eqref{eq:skew_symm_upwind_SBP_SAT-multi-block}  of the nonlinear Burgers' equation around the initial condition \(U(x) = u(x, 0)\). Subsequently, we computed the eigenvalues of the resulting linear operator \(\mathbf{Q}\). The eigenvalue spectra of \(\mathbf{Q}\) are presented in Figures \ref{fig:burgers-4th-order-eigenvalues} and \ref{fig:burgers-5th-order-eigenvalues} for various DP SBP operators, including the recently developed DP DRP FD SBP operators \cite{williams2024drp} and the DP DG SBP operators introduced in \cite{GLAUBITZ2025113841}.
To guarantee nonlinear stability, we employed the skew-symmetric form with the parameter \(\alpha = 2/3\). We conducted stability tests for different volume upwind parameters \(\gamma\), specifically the case \(\gamma = 0\) (no volume upwinding) and the case \(\gamma = \gamma_{\text{opt}} > 0\) (non-vanishing volume upwinding). The eigenvalue analyses indicate that choosing \(\gamma = \gamma_{\text{opt}} > 0\) prevents the emergence of eigenvalues with positive real parts, that is \(\Re{\lambda(\mathbf{Q})} \le 0\), which is a necessary condition for the linear stability of the entropy-stable DP SBP approximation \eqref{eq:skew_symm_upwind_SBP_SAT-multi-block} with \(\alpha = 2/3\).
Conversely, when \(\alpha = 2/3\) and the volume upwind parameter is set to zero (\(\gamma = 0\)), the eigen spectra are evenly split between the negative complex plane (\(\Re{\lambda(\mathbf{Q})} < 0\)) and the positive complex plane (\(\Re{\lambda(\mathbf{Q})} > 0\)).  Majority of eigenvalue spectra that exhibit positive real eigenvalues (\(\Re{\lambda(\mathbf{Q})} > 0\)) correspond to numerical modes with exponentially growing amplitudes. This is indicative of unstable and non-physical solutions. Such growing modes are primarily associated with unresolved high-frequency numerical oscillations, which can compromise the accuracy and stability of the numerical simulations. 

Furthermore, as illustrated in Figures \ref{fig:burgers-4th-order-eigen-modes} and \ref{fig:burgers-5th-order-eigen-modes}, we have extracted the eigenmodes associated with the fastest numerical growth rates, with initial amplitudes on the order of \(\sim 10^{-3}\). These eigenmodes serve as  perturbations of the initial conditions to  the semi-discrete approximations of the nonlinear Burgers' equation. The numerical solution is then evolved from these perturbed initial conditions using a fixed timestep of \(\Delta t = 10^{-4}\) up to the final time \(t = 10\). All time dependent experiments use the 5 stage 4th order accurate strong stability preserving RK method by Ruuth~\cite{Ruuth2006_SSPRK}.
The evolution of the perturbation norm is depicted in Figures \ref{fig:burgers-4th-order-eigen-mode-evolution} and \ref{fig:burgers-5th-order-eigen-mode-evolution}. For the case where \(\alpha = 1\) and \(\gamma = 0\), the perturbation norm remains essentially constant throughout the simulation, consistent with the linear stability analysis and specifically with the discrete estimate given by equation \eqref{eq:estimate-general-base-flow-strict-discrete}. 
In contrast, when \(\alpha = 2/3\) and \(\gamma = 0\), the perturbation norm exhibits exponential growth at a rate consistent with the eigenvalue-derived numerical growth rate. Once the amplitude of the perturbation reaches that of the initial condition, nonlinear effects become significant, leading to a saturation of the perturbation norm. At this point, the solution no longer exhibits exponential growth, and the norm stabilizes at the amplitude comparable to the initial unperturbed data. Despite this stabilization, the linear instability has effectively amplified the initial error, resulting in widespread solution contamination.
Finally, for the case \(\alpha = 2/3\) with \(\gamma = \gamma_{\text{opt}} > 0\), the perturbation norm diminishes over time, remaining small throughout the simulation duration. This behavior aligns with the eigenvalue analysis presented previously, as shown in Figures \ref{fig:burgers-4th-order-eigenvalues} and \ref{fig:burgers-5th-order-eigenvalues}, indicating linear stability under these parameters.
\begin{figure}[htbp]
    \centering
    \newcommand{\figurescale}{0.99}
    \begin{subcaptionblock}{\textwidth}
        \centering
        \includegraphics[scale=\figurescale, trim={0em 3.4em 0em 0em}, clip]{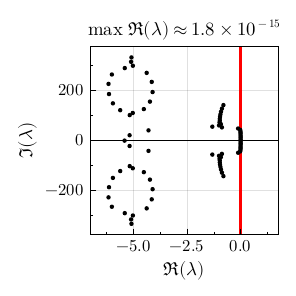}
        \includegraphics[scale=\figurescale, trim={0em 3.4em 0em 0em}, clip]{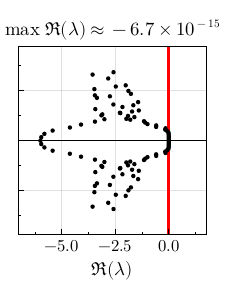}
        \includegraphics[scale=\figurescale, trim={0em 3.4em 0em 0em}, clip]{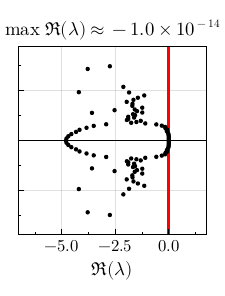}
        \\
        \includegraphics[scale=\figurescale]{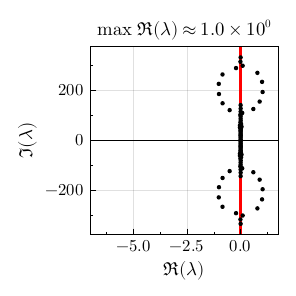}
        \includegraphics[scale=\figurescale]{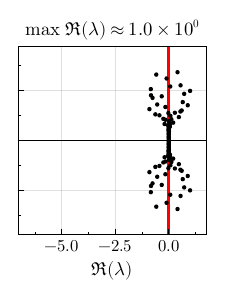}
        \includegraphics[scale=\figurescale]{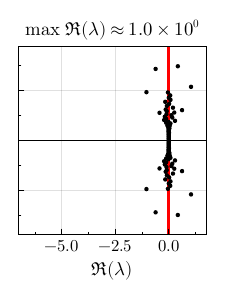}
        \caption{The eigen spectra of the linearized discrete operators. Top: with volume upwinding $\Gamma >0$. Bottom: without  volume upwinding $\Gamma =0$.}
        \label{fig:burgers-4th-order-eigenvalues}
    \end{subcaptionblock}
    \begin{subcaptionblock}{\textwidth}
        \centering
        \includegraphics[scale=\figurescale]{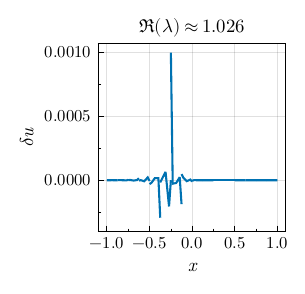}
        \includegraphics[scale=\figurescale]{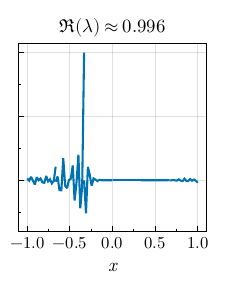}
        \includegraphics[scale=\figurescale]{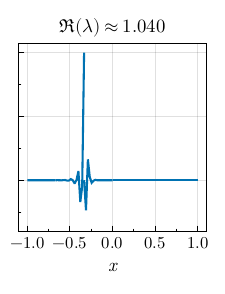}
        \caption{The eigenmode of the skew-symmetric residual.}
        \label{fig:burgers-4th-order-eigen-modes}
    \end{subcaptionblock}
    \begin{subcaptionblock}{\textwidth}
        \centering
        \includegraphics[scale=\figurescale]{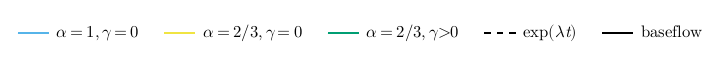} \\[-8pt]
        \includegraphics[scale=\figurescale]{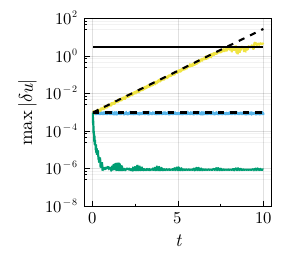}
        \includegraphics[scale=\figurescale]{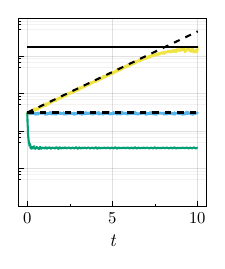}
        \includegraphics[scale=\figurescale]{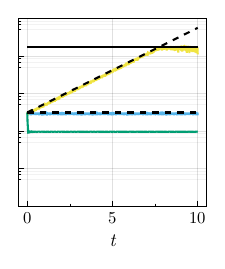}
        \caption{The evolution of the norm of  of the error over time.}
        \label{fig:burgers-4th-order-eigen-mode-evolution}
    \end{subcaptionblock}
    \caption{Some instructive numerical  experiments for the Burgers' equation investigating linear stability. Left: 4th degree polynomial DP DG approximation on 16 elements; Center: 4th order accurate DP FD approximation on 6 elements with 16 nodes each; Right: 4th order accurate DRP DP FD approximation on 6 elements with 16 nodes each.}
    \label{fig:burgers-4th-order}
\end{figure}
\begin{figure}[htbp]
    \centering
    \newcommand{\figurescale}{0.99}
    \begin{subcaptionblock}{\textwidth}
        \centering
        \includegraphics[scale=\figurescale, trim={0em 3.4em 0em 0em}, clip]{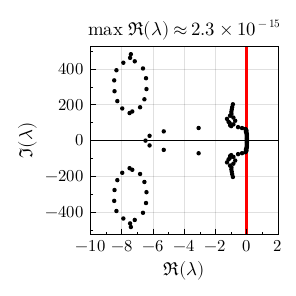}
        \includegraphics[scale=\figurescale, trim={0em 3.4em 0em 0em}, clip]{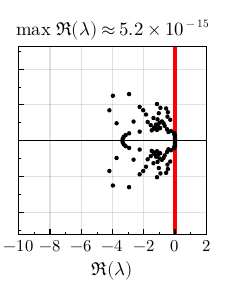}
        \includegraphics[scale=\figurescale, trim={0em 3.4em 0em 0em}, clip]{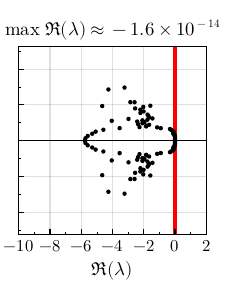}
        \\
        \includegraphics[scale=\figurescale]{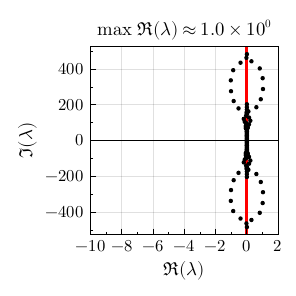}
        \includegraphics[scale=\figurescale]{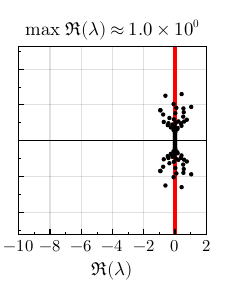}
        \includegraphics[scale=\figurescale]{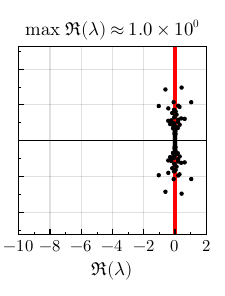}
        \caption{The eigen spectra of the linearized discrete operators. Top: with volume upwinding $\Gamma >0$. Bottom: without  volume upwinding $\Gamma =0$.}
         \label{fig:burgers-5th-order-eigenvalues}
    \end{subcaptionblock}
    \begin{subcaptionblock}{\textwidth}
        \centering
        \includegraphics[scale=\figurescale]{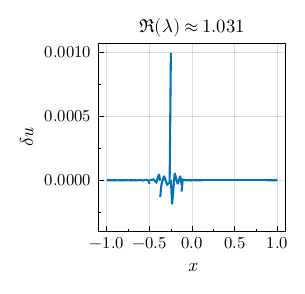}
        \includegraphics[scale=\figurescale]{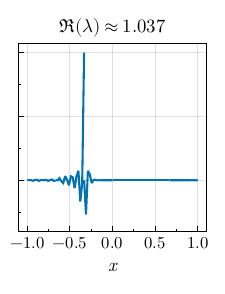}
        \includegraphics[scale=\figurescale]{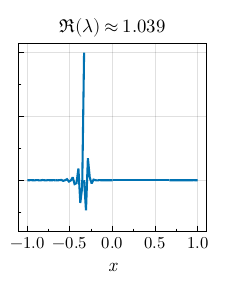}
        \caption{The eigenmode of the skew-symmetric residual.}
         \label{fig:burgers-5th-order-eigen-modes}
    \end{subcaptionblock}
    \begin{subcaptionblock}{\textwidth}
        \centering
        \includegraphics[scale=\figurescale]{plots/burgers-1D/eigenmodes/legend-horizontal.pdf} \\[-8pt]
        \includegraphics[scale=\figurescale]{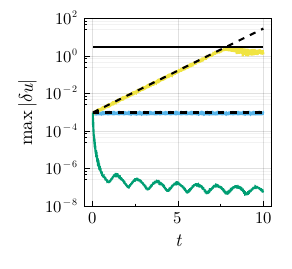}
        \includegraphics[scale=\figurescale]{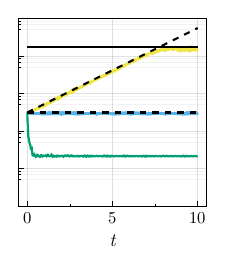}
        \includegraphics[scale=\figurescale]{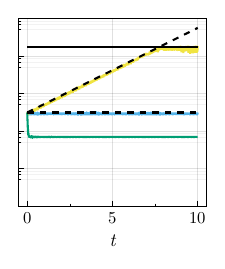}
        \caption{The evolution of the norm of  of the error over time.}
        \label{fig:burgers-5th-order-eigen-mode-evolution}
    \end{subcaptionblock}
    \caption{Some instructive numerical  experiments for the Burgers' equation investigating linear stability. Left: 5th degree polynomial DP DG approximation on 16 elements; Center: 5th order accurate DP FD approximation on 6 elements with 16 nodes each; Right: 5th order accurate DRP DP FD approximation on 6 elements with 16 nodes each.}
    \label{fig:burgers-5th-order}
\end{figure}
\section{Systems of nonlinear conservation laws: SWEs}\label{sec:swes-stability}
Next, we extend the linear stability analysis to systems of nonlinear conservation laws. Specifically, we focus on the numerical approximation of the nonlinear SWEs in both flux/conservative form and  skew-symmetric form. In 1D, the unknown vector is \(\mathbf{q} = [h, uh]^T\), where \(h > 0\) denotes the water height, \(u\) is the flow velocity, and \(g > 0\) represents the constant acceleration due to gravity. The PDE flux function is given by \(\mathbf{f}(\mathbf{q}) = [uh, hu^2 + \frac{1}{2}gh^2]^T\), and the entropy is the total mechanical energy, expressed as \(e(\mathbf{q}) = \frac{1}{2}(hu^2 + gh^2)\)>0.
As seen in the previous section, one way to ensure local linear/energy is to apply the DP FD/DG SBP operators directly to the flux, see for examples \cite{ranocha2023highorder,GLAUBITZ2025113841}. However, while the resulting method is locally linear/energy stable, detailed  numerical experiments presented in \cite{duru2024dualpairingsummationbypartsfinitedifference,ranocha2023highorder,GLAUBITZ2025113841} show that method lacks the robustness that is necessary for simulating nonlinear phenomena such as turbulence and shocks. The DP SBP entropy stable scheme derived in \cite{duru2024dualpairingsummationbypartsfinitedifference} was shown to be robust. Here we will investigate the local linear/energy stability properties of the method for systems of nonlinear SWEs.

An entropy-stable, multi-block upwind DP SBP method for the nonlinear SWEs--applied to the conserved variables--was developed in \cite{duru2024dualpairingsummationbypartsfinitedifference}. This method combines the skew-symmetric form with entropy-stable upwind splitting, and employs DP SBP operators to produce a semi-discrete approximation of the one-dimensional shallow water system that preserves key stability properties. 
The entropy stable DP FD/DG SBP approximation reads
\begin{align}\label{eq:swe1D}
    \dv{}{t} \vb q +  \vb{P}(\vb q, \vb{f}(\vb{q}), \vec{D}_x) =0,
\end{align}
%%%
where  skew-symmetric SBP approximation of the divergence of the PDE flux is given by
   \begin{align*} 
\vb{F}(\vb q, \vb{f}(\vb{q}), \vec{D}_x)=
\begin{bmatrix}
      \vec{D} ({\vec h \vec u}) \\
   \frac12 \vec{D} \left({\vec h \vec u}^2 \right) + \frac12 {\vec u} \vec{D} \left({\vec h \vec u} \right) + \frac12 {\left(\vec h \vec u\right)} \vec{D} {\vec u} + { g \vec h} \vec{D} {\vec h} 
\end{bmatrix},
\quad
\vec{g} = \begin{bmatrix}
 \vec u\\
g\vec h  -\frac12 \vec u ^2
\end{bmatrix},
\end{align*}
and spatial discrete operator with interface and volume upwinding is 
\begin{align*} 
\vb{P}(\vb q, \vb{f}(\vb{q}), \vec{D}_x) =\vb{F}(\vb q, \vb{f}(\vb{q}), \vec{D}_x) -  \frac{1}{2}\mathbf{H}_x^{-1}\left(\boldsymbol{\Upsilon} \otimes\left(\vec{B}_{Ix}+\vec{B}_{nx}\right)\right)\mathbf{g} -   \frac{1}{2}\Gamma\left(I_m \otimes\left(\vec{D}_{+}-\vec{D}_{-}\right)\right)\mathbf{g}.
\end{align*}
Here, for the $k$th element, $\Upsilon^{k}=\fn{\diag}{[\upsilon_1^k+\upsilon_1^{k+1}, \upsilon_2^k+\upsilon_2^{k+1}]}$ is the interface upwind parameter matrix and 
$\Gamma^k = \fn{\diag}{[\gamma_1^k, \gamma_2^k]}$ is the volume upwind parameter matrix with, e.g., $\upsilon_1^k = \gamma_1^k = \max_{1\le j \le n} h^k_j(|u^k| + \sqrt{gh^k_j})/(h^k_jg + (u^k_j)^2) \ge 0$, $\upsilon_2^k=\gamma_2^k=  \max_{1 \le j \le n}| h_j^ku_j^k |\ge 0$.
We define the total semi-discrete entropy $E_h(t):=\langle\mathbf{1}, \mathbf{e} \rangle_{H}
            = \mathbf{1}^T{H}\mathbf{e}$. Then it follows that 
\begin{align*}
   \dv{}{t} E_h(t)= - 
   \frac 1 2\sum_{k=1}^{K}\sum_{i=1}^2\frac{\upsilon_i^{k} + \upsilon_i^{k+1}}{2}\lJump{\mathbf{g}_i^k}\rJump^2
   -\frac 1 2  \sum_{k=1}^K\sum_{i=1}^2{\gamma_i^k\inp{ \mathbf{g}_i^k} {\left(D_{-}-D_{+} \right)\mathbf{g}_i^k}_{H}} 
      \le  0.
\end{align*}
The numerical method \eqref{eq:swe1D} is entropy stable for all $\upsilon_i^k\ge0$, $\gamma_i^k\ge 0$. We refer the reader to \cite{duru2024dualpairingsummationbypartsfinitedifference} for proofs and more elaborate discussion. It is significant to note that the skew-symmetric reformulation can support linearly unstable numerical modes which will corrupt numerical simulations everywhere. Thus, our main objective here is to investigate the local linear stability properties of the nonlinear semi-discrete approximation \eqref{eq:swe1D}.

\subsection{Local linear stability analysis for the 1D SWE}
We consider the smooth base-flow
%\begin{align*}
    $h(x)  = 8 + \sin(2\pi(x - 0.7))$, 
    $hu(x) = 1 + \sin(2\pi(x + 0.7))$, in the unit interval
    $x \in [0, 1]$, with the gravitational acceleration $g =1$,
%\end{align*}
and linearize the nonlinear semi-discrete system \eqref{eq:swe1D} around the base-flow sampled on the grid. We
obtain
\begin{align}\label{swe-linear-conservation-law-discrete}
    \dv{}{t} \vb{\delta q} = \vb{Q} \vb{\delta q}, \quad  \vb{\delta q}(0) = \vb{\delta q}_0, \quad t \in [0, T].
\end{align}
Here, $\vb{Q}$ is the Jacobian  matrix for the linearised operator, which we compute by auto-differentiating the semi-discrete approximations using the Julia package \textit{ForwardDiff.jl}~\cite{RevelsLubinPapamarkou2016}.
The eigenvalues of the discrete linear operator \(\mathbf{Q}\) are displayed in Figure \ref{fig:swe-4th-5th-order} for various DP SBP operators \cite{MATTSSON2017upwindsbp}, including the DP DRP FD SBP operators \cite{williams2024drp} and the DP DG SBP operators from \cite{GLAUBITZ2025113841}. Note that the entropy-stable skew-symmetric form with volume upwinding (\(\Gamma > 0\)) has no eigenvalues with positive real parts, which ensures local linear stability. In contrast, and similar to the 1D Bugers' equation, the eigen spectra for  the entropy-stable DGSEM/SBP approximation of the skew-symmetric form without volume upwinding (\(\Gamma = 0\))  are evenly split between the negative complex plane (\(\Re{\lambda(\mathbf{Q})} < 0\)) and the positive complex plane (\(\Re{\lambda(\mathbf{Q})} > 0\)). The eigenvalues with positive real parts indicate local linear instability. Thus, when \(\Gamma = 0\), the numerical method tends to amplify grid-scale errors, which may lead to the pollution of the numerical solutions throughout the domain.
\begin{figure}[htbp]
    \centering
    \newcommand{\figurescale}{0.9}
    \begin{subcaptionblock}{\textwidth}
        \centering
        \includegraphics[scale=\figurescale]{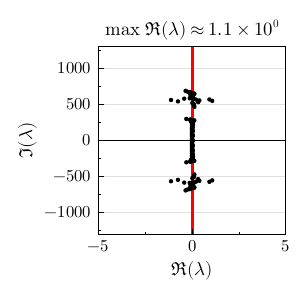}
        \includegraphics[scale=\figurescale]{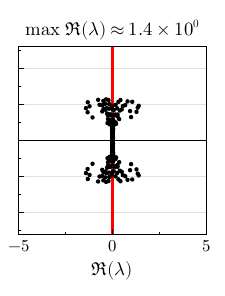}
        \includegraphics[scale=\figurescale]{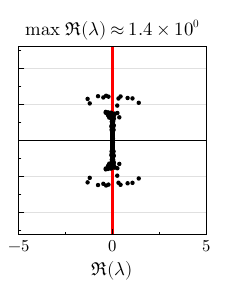}
        \\
        \includegraphics[scale=\figurescale]{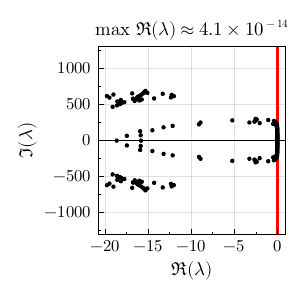}
        \includegraphics[scale=\figurescale]{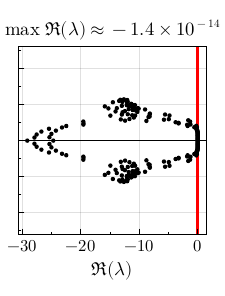}
        \includegraphics[scale=\figurescale]{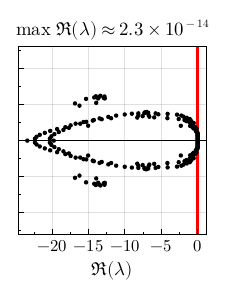}
        \caption{Left: 4th degree polynomial DP DG approximation on 16 elements. Center: 4th order accurate DP SBP on 6 elements with 16 nodes each. Right: 4th order accurate DRP DP on 6 elements with 16 nodes each. Top: without volume upwinding $\Gamma =0$. Bottom: with  volume upwinding $\Gamma >0$.}
    \end{subcaptionblock}
    \begin{subcaptionblock}{\textwidth}
        \centering
        \includegraphics[scale=\figurescale]{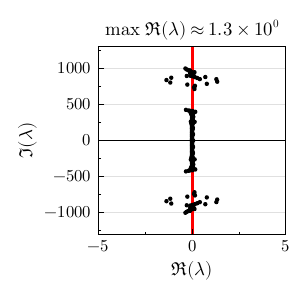}
        \includegraphics[scale=\figurescale]{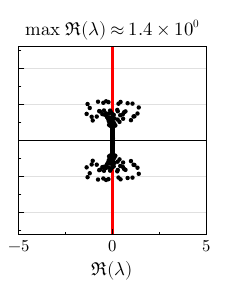}
        \includegraphics[scale=\figurescale]{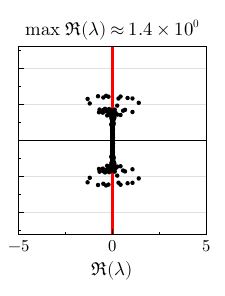}
        \\
        \includegraphics[scale=\figurescale]{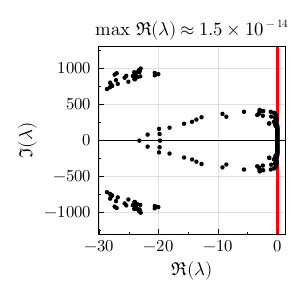}
        \includegraphics[scale=\figurescale]{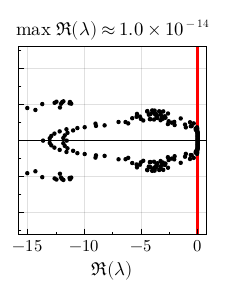}
        \includegraphics[scale=\figurescale]{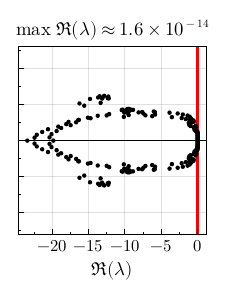}
        \caption{Left: 5th degree polynomial DP DG approximation on 16 elements. Center: 5th order accurate DP SBP on 6 elements with 16 nodes each. Right: 5th order accurate DRP DP on 6 elements with 16 nodes each. Top: without volume upwinding $\Gamma =0$. Bottom: with  volume upwinding $\Gamma >0$.}
    \end{subcaptionblock}
    \caption{
        The eigen spectra of the linearized discrete operators for the semi-discrete approximations of the  1D nonlinear SWE.
    }
    \label{fig:swe-4th-5th-order}
\end{figure}

%\subsubsection{Evolution of numerical modes}
As previously described, for different operators, we have extracted the eigenmodes for both \(\delta h\) and \(\delta hu\), corresponding to the fastest numerical growth rate, with an initial amplitude of approximately \(10^{-4}\), see Figure \ref{fig:swe-eigenmodes-4th-order}--\ref{fig:swe-eigenmodes-5th-order}. 
%%%
\begin{figure}[htbp]
     \centering
    \newcommand{\figurescale}{0.9}
    \begin{subcaptionblock}{\textwidth}
        \centering
        \includegraphics[scale=\figurescale]{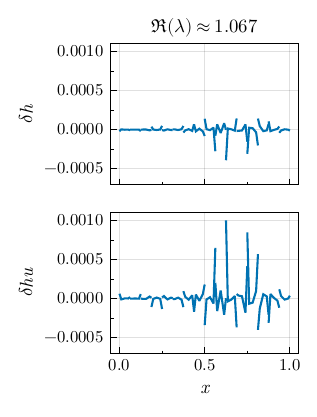}
        \includegraphics[scale=\figurescale]{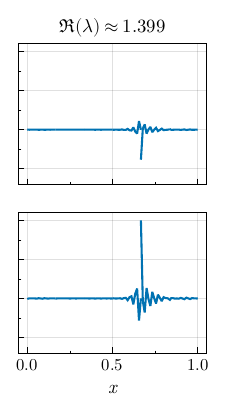}
        \includegraphics[scale=\figurescale]{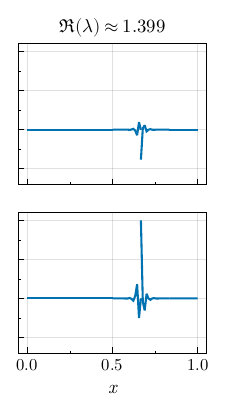}
        \caption{The eigenmode of the skew-symmetric residual.}
    \end{subcaptionblock}
    \begin{subcaptionblock}{\textwidth}
        \centering
        \includegraphics{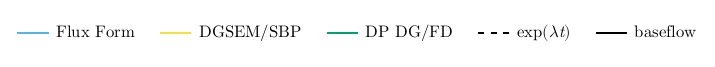} \\
        \includegraphics[scale=\figurescale]{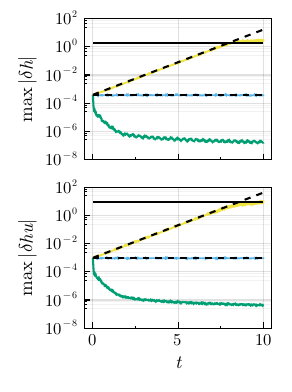}
        \includegraphics[scale=\figurescale]{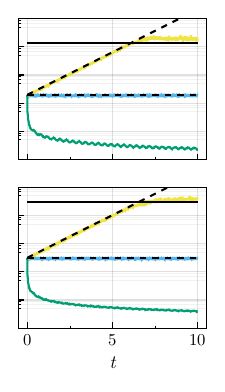}
        \includegraphics[scale=\figurescale]{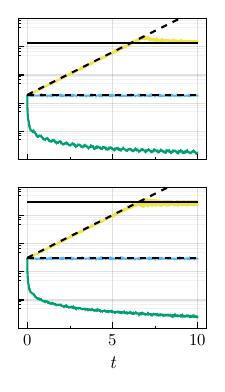}
        \caption{The evolution of the norm of  of the error over time.}
    \end{subcaptionblock}
    \caption{Numerical  experiments for the 1D SWE  investigating linear stability. Left: 4th degree polynomial DP DG approximation on 16 elements; Center: 4th order accurate DP FD approximation on 6 elements with 16 nodes each; Right: 4th order accurate DRP DP FD approximation on 6 elements with 16 nodes each.}
    \label{fig:swe-eigenmodes-4th-order}
\end{figure}
%%%
\begin{figure}[htbp]
     \centering
    \newcommand{\figurescale}{0.9}
    \begin{subcaptionblock}{\textwidth}
        \centering
        \includegraphics[scale=\figurescale]{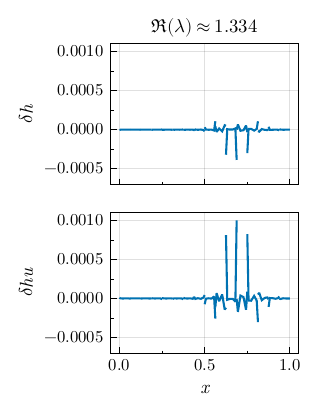}
        \includegraphics[scale=\figurescale]{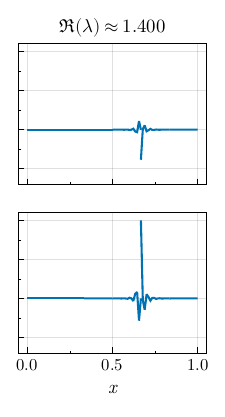}
        \includegraphics[scale=\figurescale]{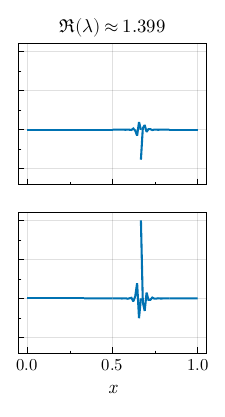}
        \caption{The eigenmode of the skew-symmetric residual.}
    \end{subcaptionblock}
    \begin{subcaptionblock}{\textwidth}
        \centering
        \includegraphics{plots/swe-1D/eigenmodes/legend-horizontal.pdf} \\
        \includegraphics[scale=\figurescale]{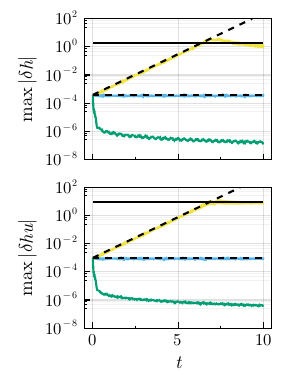}
        \includegraphics[scale=\figurescale]{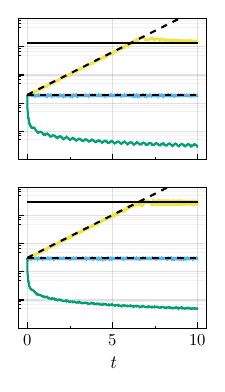}
        \includegraphics[scale=\figurescale]{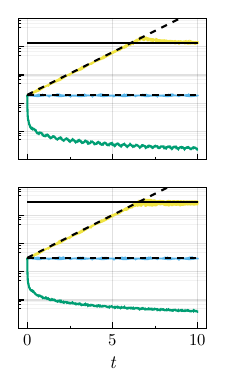}
        \caption{The evolution of the norm of  of the error over time.}
    \end{subcaptionblock}
    \caption{Numerical  experiments for the 1D SWE  investigating linear stability. Left: 5th degree polynomial DP DG approximation on 16 elements; Center: 5th order accurate DP FD approximation on 6 elements with 16 nodes each; Right: 5th order accurate DRP DP FD approximation on 6 elements with 16 nodes each.}
    \label{fig:swe-eigenmodes-5th-order}
\end{figure}
% 
%%%
We perturbed the initial conditions using these eigenmodes and evolved the numerical solutions using a fixed timestep of \(\Delta t = 10^{-5}\) until the final time \(t = 10\). The evolution of the norms of error are shown in Figure~\ref{fig:swe-eigenmodes-4th-order}--\ref{fig:swe-eigenmodes-5th-order}.
%%%
For the flux/conservative form, with \(\Gamma \geq 0\), the numerical method is provably linearly stable. Consequently, when the upwind parameter vanishes \(\Gamma = 0\), the error norm remains constant throughout the entire simulation. This observation aligns with the analysis presented in the previous section, particularly the discrete estimate \eqref{eq:estimate-general-base-flow-strict-discrete}. However, there is no proof of nonlinear or entropy stability. Therefore, it is important to note that this scheme lacks robustness and is likely to fail when encountering strongly nonlinear features such as shocks or turbulence. For more detailed results, see \cite{duru2024dualpairingsummationbypartsfinitedifference}.
% 

%%%%
When \(\Gamma = 0\), the entropy-stable skew-symmetric form becomes locally energy unstable. The norms of the errors  grow over time, with their growth rates aligning with the numerical growth rates dictated by the corresponding eigenvalues. Similarly, once the amplitude of the errors reaches that of the initial conditions, nonlinear stability is activated. At this point, the norm of the evolving numerical solutions ceases to grow and stabilizes, maintaining the amplitude of the unperturbed initial data. However, the perturbations have grown sufficiently large and have polluted the numerical simulations everywhere.
% 
%%%

%%%
Lastly, with volume upwinding parameter \(\Gamma > 0\), the skew-symmetric form is both entropy/nonlinearly stable and locally energy stable. The error norm decays over time and remains small throughout the entire simulation duration. This behavior is consistent with the numerical eigenvalue analysis presented earlier, as illustrated in  Figure \ref{fig:swe-eigenmodes-4th-order}--\ref{fig:swe-eigenmodes-5th-order}.
  % .

\subsection{Local linear stability analysis for the 2D SWE}
We extend the local linear stability analysis to the 2D SWEs. In 2D, the unknown vector field  is \(\mathbf{q} = [h, uh, vh]^T\), where \(h > 0\) denotes the water height, \([u,v]^T\) is the flow velocity vector, and \(g > 0\) represents the constant acceleration due to gravity. The PDE flux function is given by \(\mathbf{f}(\mathbf{q})=[\vb{f}_x, \vb{f}_y]\), where \(\mathbf{f}_x(\mathbf{q}) = [hu, hu^2 + \frac{1}{2}gh^2, huv]^T\) and \(\mathbf{f}_y(\mathbf{q}) = [hv, huv,  hu^2 + \frac{1}{2}gh^2]^T\). The entropy is the total mechanical energy, expressed as \(e(\mathbf{q}) = \frac{1}{2}(hu^2 + hv^2 + gh^2)\)>0. We consider the corresponding 2D skew-symmetric form  and extend the entropy stable  1D semi-discrete approximation \eqref{eq:swe1D} to 2D using tensor products. For details, please see \cite{duru2024dualpairingsummationbypartsfinitedifference}.

We consider the gravitational acceleration $g = 1$ and  linearize the  semi-discrete approximation of the 2D nonlinear SWEs around the stationary vortex given by
\begin{equation} \label{eq:stationary-vortex-2d}
    h = 12 + \int_r^\infty \frac{u^2+v^2}{sg} \dd{s}, \quad
    u = - \partial_y \psi, \quad
    v =  \partial_x \psi
\end{equation}
where
\[
    \psi(r) = \exp(-r^2), \quad r = \sqrt{x^2+y^2}, \quad (x,y)\in [-\pi, \pi]^2.
\]
As above the corresponding  linear operator \(\vb{Q}\) of the nonlinear semi-discrete approximation is computed using auto-differentiation. 
The eigenvalues of the discrete linear operator \(\mathbf{Q}\) are illustrated in Figure \ref{fig:swe-2D-5th-order} for various DP SBP operators, including the DP DRP FD SBP operators \cite{williams2024drp} and the DP DG SBP operators from \cite{GLAUBITZ2025113841}.

For the entropy-stable, skew-symmetric form without volume upwinding (\(\Gamma = 0\)), the eigenvalue spectra are symmetrically distributed between the negative and positive halves of the complex plane (\(\Re{\lambda(\mathbf{Q})} < 0\) and \(\Re{\lambda(\mathbf{Q})} > 0\), respectively). Most eigenvalues with positive real parts (\(\Re{\lambda(\mathbf{Q})} > 0\)) correspond to numerical modes exhibiting exponential growth. These eigenvalues, with \(\Re{\lambda} \sim 4 \times 10^{-1}\), indicate the presence of unstable numerical modes. Such growing modes are primarily associated with unresolved high-frequency oscillations, which can undermine the accuracy and stability of the simulations. To substantiate this, we plotted the wave mode associated with the largest real eigenvalues for the DGSEM in Figure \ref{fig:swe-2D-stationary-vortex-eigenmodes}. As shown, this mode is unresolved and indicative of numerical instability. When \(\Gamma = 0\), the numerical scheme tends to amplify grid-scale errors, leading to deterioration in solution quality across the domain and preventing convergence.

In contrast, the entropy-stable DP DG scheme with volume upwinding (\(\Gamma > 0\)) demonstrates local linear stability, as most eigenvalues have non-positive real parts (\(\Re{\lambda(\mathbf{Q})} \le 0\)). However, some eigenvalues still possess positive real parts, with the largest \(\Re{\lambda} \sim 4 \times 10^{-2}\). These modes are physical linear modes with bounded exponential growth. To verify this, we plotted the wave mode corresponding to the largest real eigenvalue for the DP DG scheme with volume upwinding in Figure \ref{fig:swe-2D-stationary-vortex-eigenmodes}. The mode appears sufficiently smooth and corresponds to a resolved physical mode with bounded growth.

Further, we perturbed the initial conditions with the unstable numerical modes exhibiting the fastest growth rates, with initial perturbation amplitudes around \(10^{-3}\). The solutions were evolved until \(t = 100\) using a fixed timestep of \(5\times10^{-3} \Delta x\). Snapshots of the numerical solutions for DGSEM and DP DG schemes are shown in Figure \ref{fig:swe-2D-stationary-vortex-snapshot}, while the evolution of the error norms is presented in Figure \ref{fig:swe-2D-stationary-vortex-eigenmodes}.

When \(\Gamma = 0\), the entropy-stable skew-symmetric form becomes locally energy unstable, causing the error norms to grow at rates aligned with the eigenvalues. Once errors reach the initial perturbation amplitude, nonlinear stability mechanisms activate, stabilizing the error norm and halting further growth. Nonetheless, this linear instability amplifies the initial perturbation, leading to widespread solution contamination, as clearly visible in the left panel of Figure \ref{fig:swe-2D-stationary-vortex-snapshot}.

Conversely, for \(\Gamma > 0\), the skew-symmetric form exhibits both entropy and nonlinear stability, with the error norm exhibits marginal growth and remains small throughout the simulation, consistent with the eigenvalue analysis in Figure \ref{fig:swe-2D-stationary-vortex-eigenmodes}. As shown in the right panel of Figure \ref{fig:swe-2D-stationary-vortex-snapshot}, the DP DG scheme with volume upwinding (\(\Gamma >0\)) accurately approximates the stationary vortex throughout the entire simulation duration.
%%0_elems=_10, 10_nodes=_6, 6_operator=GlaubitzEtal2024_-0.1_order=5
%%0_elems=_10, 10_nodes=_6, 6_operator=GlaubitzEtal2024_-0.1_order=5
%%0_elems=_10, 10_nodes=_6, 6)_operator=GlaubitzEtal2024_-0.1_order=5_scheme=Skew Symm._γ=0_
%%0_elems=_3, 3_nodes=_20, 20_operator=WilliamsDuru2024_order=5_scheme=Skew Symm_γ>0_
\begin{figure}[htbp]
    \centering
    \newcommand{\figurescale}{0.9}
    \includegraphics[scale=\figurescale]{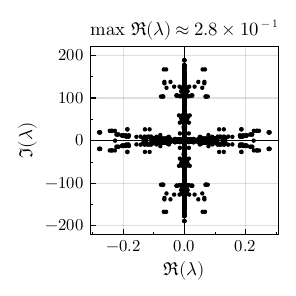}
    \includegraphics[scale=\figurescale]{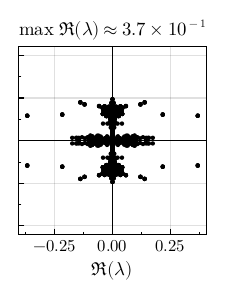}
    \includegraphics[scale=\figurescale]{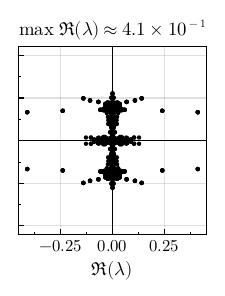}
    \\
    \includegraphics[scale=\figurescale]{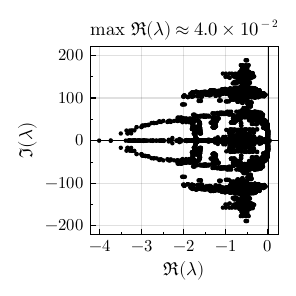}
    \includegraphics[scale=\figurescale]{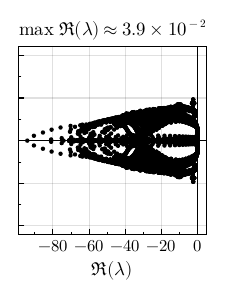}
    \includegraphics[scale=\figurescale]{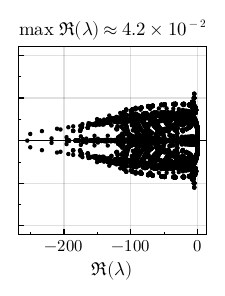}
    \caption{
    The eigen spectra of the linearised discrete operators for the semi-discrete approximations of the  2D nonlinear SWE. Top: Without volume upwinding $\Gamma =0$. Bottom: With volume  upwinding $\Gamma >0$.
    Left: 5th degree polynomials DP DG approximation on \(10^2\) elements. Center: 5th order accurate DP SBP approximations on \(3^2\) elements with \(20^2\) nodes each. Right: 5th order DRP DP on \(3^2\) elements with \(20^2\) nodes each.
    }
    \label{fig:swe-2D-5th-order}
\end{figure}
\begin{figure}
    \centering
    \newcommand{\figurescale}{0.99}
    \begin{subcaptionblock}{\textwidth}
        \centering
        \includegraphics[width=0.5\textwidth]{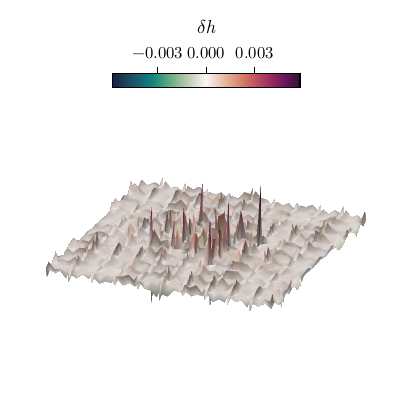}%
        \includegraphics[width=0.5\textwidth]{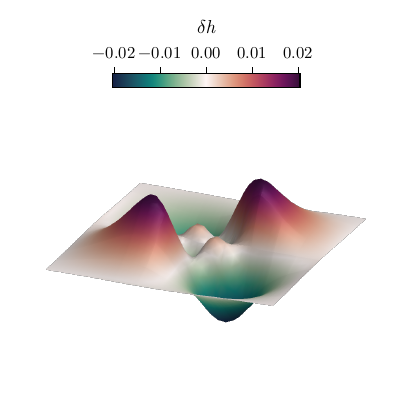}
        \caption{The fastest growing eigenmodes. Left:  The DGSEM without volume upwinding $\Gamma =0$ and eigenvalue \(\lambda=0.279 \pm 19.4i\). Right: The DP DG with  volume upwinding $\Gamma >0$ and eigenvalue \(\lambda=0.0400 \pm 0.543i\).}
    \end{subcaptionblock}
    \begin{subcaptionblock}{\textwidth}
        \centering
        \includegraphics[scale=\figurescale]{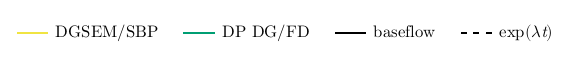}\\
        \includegraphics[scale=\figurescale]{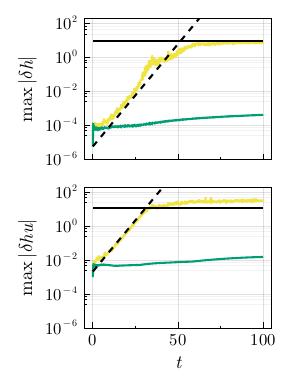}
        \includegraphics[scale=\figurescale]{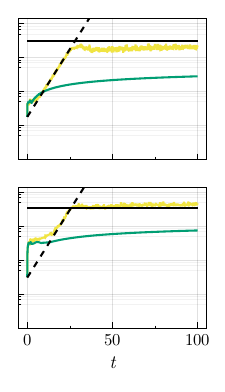}
        \includegraphics[scale=\figurescale]{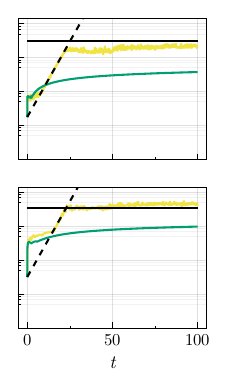}
        \caption{
            The evolution of the norm of  of the error over time.
            Left: 5th degree polynomial DP DG approximation on \(10^2\) elements; Center: 5th order accurate DP FD approximation on \(3^2\) elements with \(20^2\) nodes each; Right: 5th order accurate DRP DP FD approximation on \(3^2\) elements with \(20^2\) nodes each
        }
    \end{subcaptionblock}
    \caption{Numerical experiments for the 2D SWE stationary vortex investigating linear stability.}
    \label{fig:swe-2D-stationary-vortex-eigenmodes}
\end{figure}
\begin{figure}
    \centering
    \includegraphics{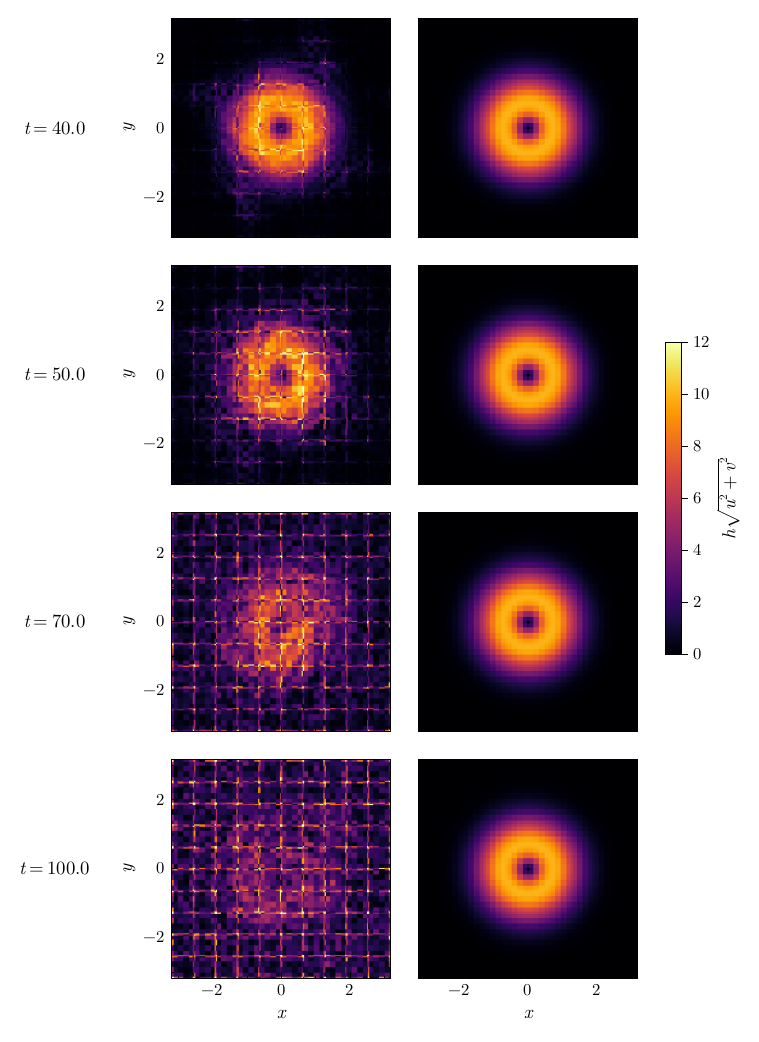}
    \caption{Stationary vortex (SWE): The snapshots of the absolute momentum for the DGSEM (left) and the DP DG method (right) using 6th degree polynomials on a \(10^2\) element mesh.}
    \label{fig:swe-2D-stationary-vortex-snapshot}
\end{figure}
\subsection{Barotropic shear instability with fully developed turbulence}
Finally, we assess the accuracy of the DP DG/FD SBP numerical framework in capturing turbulent scales. Specifically, we examine the barotropic shear instability \cite{hew2024stronglystabledualpairingsummation, galewsky2004initial, peixoto2019semi}, characterized by fully developed turbulence modeled through the 2D nonlinear SWEs. This scenario, also known as the Kelvin-Helmholtz instability, involves initiating a barotropic shear within zonal jets by introducing a sharp fluid discontinuity, augmented with small Gaussian perturbations.

The computational domain is defined as \(\Omega = [0, L]^2\), where \(L = 2\pi \times \SI{6371.22}{\km}\), with periodic boundary conditions in both directions. The initial conditions are given by:
\begin{align*}
u_0 &= \bar{u}_0 \left(\operatorname{sech}\frac{y - y_+}{10^6} - \operatorname{sech}\frac{y - y_-}{10^6}\right), \quad
v_0 = 0, \quad
h_0 = H - \frac{f}{g} \int_0^y u(x, s) \, ds + \tilde{h}(x, y),
\end{align*}
where
\begin{align*}\tilde{h}(x, y) = \bar{h}_0 \sum_{i=1}^2 \exp\left(-k\, d_i(x, y)\right), \quad
d_i(x, y) = \frac{(x - x_i)^2}{L^2} + \frac{(y - y_i)^2}{L^2},\end{align*}
with:
\(\bar{h}_0 = 0.01 H\),
\((x_1, y_1) = (0.15 L, y_+)\),
\((x_2, y_2) = (0.85 L, y_-)\),
\(y_+ = 0.25 L\),
\(y_- = 0.75 L\),
\(\bar{u}_0 = \SI{50}{\m\per\s}\),
 \(f = 7.292 \times 10^{-5}\),
 \(g = \SI{9.80616}{\m\per\square\s}\),
 \(H = \SI{10}{\km}\),
and \(k = 10^3\).

We employ the entropy-stable DGSEM and the DP DG SBP scheme \cite{GLAUBITZ2025113841}, both using sixth-degree polynomials on a \(64^2\) element mesh. Similar numerical results have been obtained with DP FD SBP operators \cite{williams2024drp, MATTSSON2017upwindsbp}, although these are not shown here. The simulations are advanced with a fixed timestep \(\Delta t = 2 \times 10^{-3} \Delta x\) until the final time \(t=80\,\text{days}\). Snapshots of the absolute vorticity at various times are displayed in Figure \ref{fig:swe-2D-kh-snapshots}.
\begin{figure}[htbp]    
    \centering    
    \includegraphics{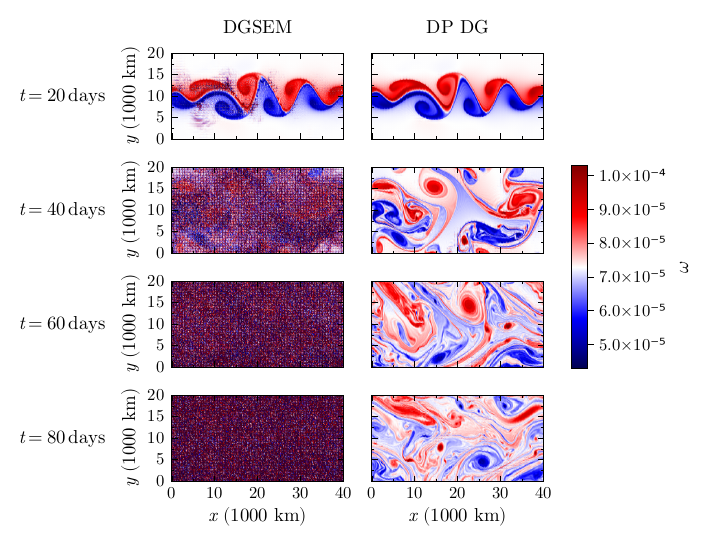}    
    \caption{Kelvin-Helmholtz instability (SWE): The snapshots of the absolute vorticity for the DGSEM (left) and the DP DG method (right) using 6th degree polynomials on a \(64^2\) element mesh.}    
    \label{fig:swe-2D-kh-snapshots}
\end{figure}
The entropy-stable DGSEM without volume upwinding (\(\Gamma=0\)) introduces numerical noise that rapidly pollutes the solution across the domain. In contrast, the entropy-stable DP DG method with volume upwinding (\(\Gamma > 0\)) effectively captures the vortex structures without generating spurious high-frequency oscillations. These observations are consistent with the theoretical analysis and the numerical results discussed in the preceding sections.

\begin{figure}[htbp]
    \centering
    \begin{subcaptionblock}{\textwidth}
        \centering
        \includegraphics{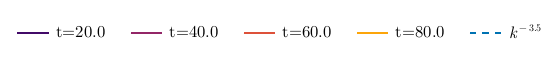}
        \\
        \includegraphics{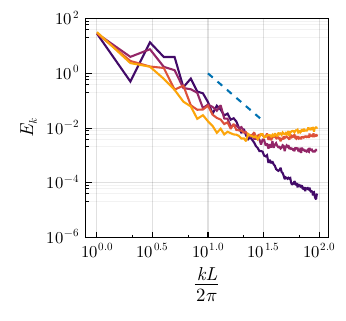}
        \includegraphics{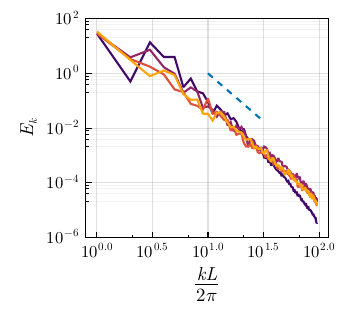}
    \end{subcaptionblock}
    \begin{subcaptionblock}{\textwidth}
        \centering
        \includegraphics{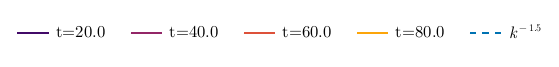}
        \\
        \includegraphics{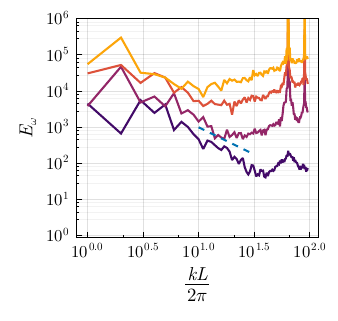}
        \includegraphics{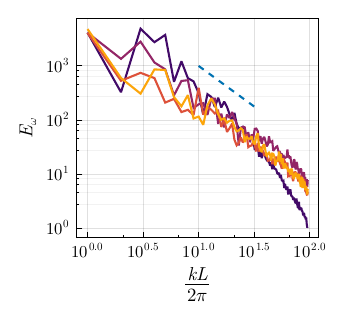}
    \end{subcaptionblock}
    \caption{The kinetic energy spectra (top) and the enstrophy spectra (bottom) for the DGSEM without volume upwinding (left) and the DP DG method with volume upwinding (right) using 6th degree polynomials on a \(64^2\) element mesh.}
    \label{fig:swe-2D-kh-spectra-dg}
\end{figure}

To demonstrate the effectiveness of the DP-DG framework in resolving turbulent scales, we analyze the shell-integrated kinetic energy and enstrophy spectra in our simulations (see Figure~\ref{fig:swe-2D-kh-snapshots}). The kinetic energy spectrum is defined as:
\begin{align*}    E_n = \sum_{n \leq |\mathbf{k}| < n + 1} E_{\mathbf{k}}, \quad \text{where} \quad E_{\mathbf{k}} = \frac{1}{2} \left( \hat{u}(\mathbf{k}) \hat{u}^*(\mathbf{k}) + \hat{v}(\mathbf{k}) \hat{v}^*(\mathbf{k}) \right),\end{align*}
with $(\hat{u}(\mathbf{k}), \hat{v}(\mathbf{k}))$ representing the Fourier-transformed velocity components, $\mathbf{k} = (k_1, k_2)$ denoting a wavenumber vector, and $^*$ indicating complex conjugation. Here, $|\mathbf{k}| = \sqrt{k_1^2 + k_2^2}$ is the magnitude of the wavenumber, and $E_n$ reflects the spectral density corresponding to wavenumber $n$, which relates to a wavelength of $L/n$, with $L$ being the domain size. 
The enstrophy spectrum, assuming isotropy, is defined as $E_\omega(k) = k^2 E_k$, consistent with the expected power-law cascades observed in two-dimensional hydrodynamic turbulence \cite{boffetta2012two}.

For the DGSEM without volume upwinding ($\Gamma = 0$), after 20 days, the simulation is completely overwhelmed by numerical noise. Consequently, the numerical spectra of kinetic energy and enstrophy do not exhibit any physically meaningful trend or power-law behavior. 
In contrast, for the DP DG scheme with volume upwinding ($\Gamma > 0$), the kinetic energy spectra--shown in the top panel of Figure ~\ref{fig:swe-2D-kh-spectra-dg}--demonstrate a reasonable Kraichnan inertial range characterized by a $k^{-3.5}$ power law \cite{Kraichnan:1968}, up to scales affected by small-scale intermittency or finite resolution effects. This indicates that the DP DG scheme effectively resolves the direct kinetic energy cascade despite having somewhat lower resolution compared to most canonical 2D turbulence studies (e.g., \cite{boffetta2010evidence}).
Similarly, the enstrophy spectra (Figure ~\ref{fig:swe-2D-kh-spectra-dg}) display a reasonable power-law behavior of approximately $k^{-1.5}$ within a narrow inertial range. It is important to note that accurately capturing the enstrophy cascade--being a second moment of velocity--generally requires high spatial resolution to achieve good convergence within the inertial range.

\section{Conclusions}\label{sec:conclusion}
Three desirable properties of numerical methods for nonlinear conservation laws are \cite{ranocha2023highorder}:
1) High-order accuracy,
2) Nonlinear (entropy) stability,
3) Local linear-energy stability.
An open question remains whether higher-order methods can simultaneously possess  these two stability properties. Nonlinear or entropy stability enhances robustness, allowing simulations to run without crashing \cite{duru2024dualpairingsummationbypartsfinitedifference}. However, because of the amplification of grid-scale errors from unresolved nonlinear wave modes, numerical simulations  might though be stable but still remain inaccurate or fail to converge, even for smooth solutions.  Meanwhile, local linear-energy stability ensures that grid-scale errors do not dominate the numerical solution.
As demonstrated in \cite{gassner2022stability}, many split-form, entropy-stable high-order schemes fail to achieve local linear-energy stability and thus lack this critical property. 

In this work, we have investigated the local linear stability of the DP DG/FD SBP schemes recently introduce in \cite{duru2024dualpairingsummationbypartsfinitedifference} within the split-form framework. Our analysis suggests that the entropy-stable volume upwind filter inherent in these methods can guarantee local linear stability. The DP DG/FD SBP framework introduces a novel numerical approach for designing reliable high-order methods for nonlinear conservation laws that are both provably entropy-stable and locally energy-stable. These theoretical insights are supported by numerical experiments involving the inviscid Burgers’ equation and nonlinear SWEs in both 1D and 2D. Furthermore, we present accurate simulations of 2D barotropic shear instability with fully developed turbulence, demonstrating the effectiveness of the DP SBP method in capturing turbulent scales efficiently.
\appendix

\bibliography{papers}

\bibliographystyle{plain}

\end{document}